\documentclass[12pt]{article}
\usepackage{amsmath}
\usepackage{amssymb}

\def\q{\quad}
\def\qq{\qquad}

\def\mod#1{\ (\text{\rm mod}\ #1)}
\def\t{\text}
\def\f{\frac}
\def\e{\equiv}
\def\b{\binom}

\def\sls#1#2{(\f{#1}{#2})}
 \def\ls#1#2{\big(\f{#1}{#2}\big)}
\def\Ls#1#2{\Big(\f{#1}{#2}\Big)}
\def\v2{\vskip0.2cm}

\begin{document}
 \centerline {\bf
     Supercongruences via Beukers' method}

\par\q\newline
\centerline{Zhi-Hong Sun}\newline \centerline{School of Mathematics
and Statistics}
\centerline{Huaiyin Normal University}
\centerline{Huaian, Jiangsu 223300, P.R. China} \centerline{Email:
zhsun@hytc.edu.cn} \centerline{URL:
http://maths.hytc.edu.cn/szh1.htm}
\par\q
\newline\centerline{Dongxi Ye}
\newline \centerline{School of Mathematics
(Zhuhai)} \centerline{Sun Yat-sen University}
 \centerline{Zhuhai, Guangdong 519082, P.R. China}
\centerline{Email: yedx3@mail.sysu.edu.cn}
\par\q\newline \centerline{}

\par {\it Abstract.} Recently, using modular forms F. Beukers posed a unified method that can deal with a large number of supercongruences involving binomial coefficients and Ap\'ery-like numbers. In this paper, we use Beukers' method to prove some conjectures of the first author concerning the congruences for
$$\sum_{k=0}^{(p-1)/2}\frac{\binom{2k}k^3}{m^k},
\ \sum_{k=0}^{p-1}\frac{\binom{2k}k^2\binom{4k}{2k}}{m^k},
\ \sum_{k=0}^{p-1}\frac{\binom{2k}k\binom{3k}k\binom{6k}{3k}}{m^k},
\  \sum_{n=0}^{p-1}\frac{V_n}{m^n},\ \sum_{n=0}^{p-1}\frac{T_n}{m^n},\  \sum_{n=0}^{p-1}\frac{D_n}{m^n}
$$ and $\sum_{n=0}^{p-1}(-1)^nA_n$ modulo $p^3$, where $p$ is an odd prime representable by some suitable binary quadratic form, $m$ is an integer not divisible by $p$, $V_n=\sum_{k=0}^n\binom{2k}k^2\binom{2n-2k}{n-k}^2$, $T_n=\sum_{k=0}^n\binom nk^2\binom{2k}n^2$, $D_n=\sum_{k=0}^n\binom nk^2\binom{2k}k\binom{2n-2k}{n-k}$ and $A_n$ is the Ap\'ery number given by $A_n=\sum_{k=0}^n\binom nk^2\binom{n+k}k^2$.

 \par\q

\section*{1. Introduction}

\par\q In 1998, Ono[14] established congruences
for $\sum_{k=0}^{p-1}\b{2k}k^3\f 1{m^k}$ modulo $p$ in the cases
$m=1,-8,16,-64,256,-512,4096$ for any prime $p\ne2,7$. Following these, for such values of $m$, in [25] Z.W. Sun hypothesized congruences for
$\sum_{k=0}^{p-1}\b{2k}k^3\f 1{m^k}$ modulo $p^2$, which have been sequentially
proved by the first author in [17] and Kibelbek et al in [10]. Going beyond that, in [21],
the first author also conjecturally formulated congruences for
$\sum_{k=0}^{p-1}\b{2k}k^3\f 1{m^k}$ modulo $p^3$. Before instantiating these congruences, hereafter,
for positive integers $a,b$ and $n$, if
$n=ax^2+by^2$ for some integers $x$ and $y$, we simply write that
$n=ax^2+by^2$, and recall the well known results (see [7,24]) that for an odd prime $p$,
\begin{align*} & p=x^2+4y^2\q\t{for}\q p\e 1\mod 4,
\\&p=x^2+3y^2\q\t{for}\q p\e 1\mod 3,
\\&p=x^2+2y^2\q\t{for}\q  p\e 1,3\mod 8,
\\& p=x^2+7y^2\q\t{for}\q p\e  1,2,4\mod 7,
\end{align*}
where the integers $x$ and $y$ are uniquely determined up to sign. Following these, an example of the first author's conjectures can be stated as that for
any prime $p\not=2,7$,

$$\sum_{k=0}^{p-1}\b{2k}k^3\e
\begin{cases}  4x^2-2p-\f{p^2}{4x^2}\mod {p^3}\\\qq\qq\q\t{if $p\e 1,2,4\mod 7$ and
so $p=x^2+7y^2$,}
\\
-11p^2\b{[3p/7]}{[p/7]}^{-2} \mod {p^3}\; \q\qq\qq \t{if $p\e 3\mod 7$,}
\\
-\f{11}{16}p^2\b{[3p/7]}{[p/7]}^{-2}\mod {p^3}\q\ \qq\qq \t{if $p\e 5\mod
7$,}
\\
-\f{11}4p^2\b{[3p/7]}{[p/7]}^{-2}\mod {p^3}\q\ \qq\qq\t{if $p\e 6\mod 7$,}
\end{cases}\eqno{(1.1)}$$
where $[a]$ is the greatest integer not exceeding $a$.

\par  Based on the work of Long and
Ramakrishna[9], the first author[22] illustrated that for any odd prime $p$,
$$\sum_{k=0}^{\f{p-1}2}\f{\b{2k}k^3}{64^k}\e
\begin{cases} 4x^2-2p-\f{p^2}{4x^2}\mod{p^3}& \t{if $p=x^2+4y^2\e
1\mod 4$,}\\
-p^2\b{\f{p-1}2}{\f{p-3}4}^{-2}\mod{p^3}& \t{if $p\e 3\mod
4$}\end{cases}\eqno{(1.2)}$$
and
$$(-1)^{\f{p-1}4}\sum_{k=0}^{\f{p-1}2}\f{\b{2k}k^3}{(-512)^k}\e
 4x^2-2p-\f{p^2}{4x^2}\mod{p^3}\  \t{for $p=4k+1=x^2+4y^2$}.\eqno{(1.3)}$$
\par What's more,
in [21] and [22], the first author posed many conjectural congruences for
$$\sum_{k=0}^{p-1}\f{\b{2k}k^2\b{3k}k}{m^k},\q \sum_{k=0}^{p-1}\f{\b{2k}k^2\b{4k}{2k}}{m^k}
\q\t{and}\q\sum_{k=0}^{p-1}\f{\b{2k}k\b{3k}k\b{6k}{3k}}{m^k}$$
modulo $p^3$ for prime $p>3$ and integer $m$ coprime to~$p$. See also [16-19] and [28]  for relevant congruences modulo $p^2$. For example, the first author[21] experimentally found that
$$\sum_{k=0}^{p-1}\f{\b{2k}k^2\b{4k}{2k}}{256^k}\e
\begin{cases} 4x^2-2p-\f{p^2}{4x^2}\mod{p^3}&\t{if $p=x^2+2y^2\e 1,3\mod
8$,}
\\\f{p^2}3\b{[p/4]}{[p/8]}^{-2}\mod{p^3}&\t{if $p\e 5\mod 8$,}
\\-\f 32{p^2}\b{[p/4]}{[p/8]}^{-2}\mod{p^3}&\t{if $p\e 7\mod 8$,}
\end{cases}\eqno{(1.4)}$$
$$\Ls {-3}p\sum_{k=0}^{p-1}\f{\b{2k}k\b{3k}k\b{6k}{3k}}{12^{3k}}\e
\begin{cases} 4x^2-2p-\f{p^2}{4x^2}\mod{p^3}\q \t{if $p=x^2+4y^2\e 1\mod
4$,}
\\\f 5{12}p^2\b{(p-3)/2}{(p-3)/4}^{-2}\mod{p^3}\qq\qq\qq\t{if $p\e 3\mod 4$,}
\end{cases}\eqno{(1.5)}$$
where $\ls ap$ is the Jacobi symbol. It is worth remarking that J.-C. Liu[8] formulated the supercongruences for the sums in (1.4)-(1.5) modulo $p^3$ in terms of $p$-adic Gamma functions.
In [11] , using $p$-adic Gamma functions and Jacobi sums Mao proved (1.4) in the cases $p\e 1,5,7\mod 8$. The congruences (1.4) and (1.5) modulo $p^2$ were conjectured by Rodriguez-Villegas[15] and proved by
Mortenson[13] and Z.W. Sun[26].
\par Recently, in order to conceptually interpret all these conjectures,  using modular forms Beukers[2]  found a unified way to deal with supercongruences modulo $p^2$ or $p^3$ for sums involving binomial coefficients and Ap\'ery-like numbers, where $p$ is an odd prime representable by  some suitable binary quadratic form. In this paper, using Weber functions and Beukers' method we prove a number of congruences conjectured by the first author. In particular, we rigorously verify a number of CM values of modular functions that were stated without proofs in Beukers' work [2].
\par Thanks to Beukers for his summary [2, Appendix C] of products of three binomial coefficients, Ap\'ery-like numbers and their associated Hauptmoduls, as well as findings of their related CM values, we notice that a number of conjectures of the first author can also be charted similarly. In what follows, we state all that we have been able to attain based on Beukers' work.

The following congruences that we shall prove in Sections 3 and 4 were conjectured by the first author in [21] and [22].
\vskip0.2cm
\par{\bf Theorem 1.1} {\it Let $p$ be an odd prime, $p\e 1,2,4\mod 7$ and so $p=x^2+7y^2$. Then
$$\sum_{k=0}^{(p-1)/2}\b{2k}k^3\e(-1)^{\f{p-1}2}\sum_{k=0}^{(p-1)/2}
\f{\b{2k}k^3}{4096^k}\e
4x^2-2p-\f{p^2}{4x^2}\mod {p^3}.$$}

\par{\bf Theorem 1.2} {\it Let $p$ be a prime of the form $3k+1$ and so $p=x^2+3y^2$. Then
$$\sum_{k=0}^{(p-1)/2}\f{\b{2k}k^3}{16^k}\e (-1)^{\f{p-1}2}\sum_{k=0}^{(p-1)/2}\f{\b{2k}k^3}{256^k}\e 4x^2-2p-\f{p^2}{4x^2}\mod {p^3}.$$}

\par{\bf Theorem 1.3} {\it Let $p$ be a prime of the form $4k+1$ and so $p=x^2+4y^2$. Then
$$\sum_{k=0}^{(p-1)/2}\f{\b{2k}k^3}{(-8)^k}\e  4x^2-2p-\f{p^2}{4x^2}\mod {p^3}.$$}

\par{\bf Theorem 1.4} {\it Let $p$ be a prime such that $p\e 1,3\mod 8$ and so $p=x^2+2y^2$. Then
$$\sum_{k=0}^{(p-1)/2}\f{\b{2k}k^3}{(-64)^k} \e (-1)^{\f{p-1}2}\Big(4x^2-2p-\f{p^2}{4x^2}\Big)\mod {p^3}.$$}

\par{\bf Theorem 1.5} {\it Let $p$ be a prime such that $p\e 1,3\mod 8$ and so $p=x^2+2y^2$. Then
$$\sum_{k=0}^{p-1}\f{\b{2k}k^2\b{4k}{2k}}{256^k}
\e 4x^2-2p-\f{p^2}{4x^2}\mod {p^3}.$$}

\par{\bf Theorem 1.6} {\it Let $p$ be a prime such that $p\e 1\mod 3$ and so $p=x^2+3y^2$. Then
$$\sum_{k=0}^{p-1}\f{\b{2k}k^2\b{4k}{2k}}{(-144)^k}
\e 4x^2-2p-\f{p^2}{4x^2}\mod {p^3}.$$}

\par{\bf Theorem 1.7} {\it Let $p$ be a prime such that $p\e 1\mod 4$ and so $p=x^2+4y^2$. Then
$$\sum_{k=0}^{p-1}\f{\b{2k}k^2\b{4k}{2k}}{648^k}
\e 4x^2-2p-\f{p^2}{4x^2}\mod {p^3}.$$}

\par{\bf Theorem 1.8} {\it Let $p$ be an odd prime such that $p\e 1,2,4\mod 7$ and so $p=x^2+7y^2$. Then
$$\sum_{k=0}^{p-1}\f{\b{2k}k^2\b{4k}{2k}}{81^k}\e \sum_{k=0}^{p-1}\f{\b{2k}k^2\b{4k}{2k}}{(-3969)^k}
\e 4x^2-2p-\f{p^2}{4x^2}\mod {p^3}.$$}

\par{\bf Theorem 1.9} {\it Let $p$ be a prime such that $p\e 1,3\mod 8$ and so $p=x^2+2y^2$. Then
$$\sum_{k=0}^{p-1}\f{\b{2k}k^2\b{4k}{2k}}{28^{4k}}
\e 4x^2-2p-\f{p^2}{4x^2}\mod {p^3}.$$}

\par\textbf{Theorem 1.10} {\it Let $p$ be a prime of the form $4k+1$. Then
\begin{align*}&\sum_{k=0}^{p-1}\f{\b{2k}k^2\b{4k}{2k}}{(-12288)^k}
\\&\e \begin{cases}4x^2-2p-\f{p^2}{4x^2}\mod {p^3}&\t{if $p\e 1\mod {12}$ and so $p=x^2+9y^2$,}
\\-2x^2+2p+\f{p^2}{2x^2}\mod {p^3}&\t{if $p\e 5\mod {12}$ and so $2p=x^2+9y^2$}\end{cases}\end{align*}
and
\begin{align*}&\sum_{k=0}^{p-1}\f{\b{2k}k^2\b{4k}{2k}}{(-6635520)^k}
\\&\e \begin{cases}4x^2-2p-\f{p^2}{4x^2}\mod {p^3}&\t{if $p\e 1,9\mod {20}$ and so $p=x^2+25y^2$,}
\\-2x^2+2p+\f{p^2}{2x^2}\mod {p^3}&\t{if $p\e 13,17\mod {20}$ and so $2p=x^2+25y^2$.}\end{cases}\end{align*}}

\par\textbf{Theorem 1.11} {\sl Let $m\in\{5,13,37\}$ and
$D(m)=-1024,-82944,-14112^2$ according as $m=5,13,37$. Suppose that $p$ is an odd prime such that $\sls{-m}p=1$. Then

\begin{align*}&\sum_{k=0}^{p-1}\f{\b{2k}k^2\b{4k}{2k}}{D(m)^k}
\\&\e \begin{cases}4x^2-2p-\f{p^2}{4x^2}\mod {p^3}&\t{if $\sls{-1}p=\sls mp=1$ and so $p=x^2+my^2$,}
\\-2x^2+2p+\f{p^2}{2x^2}\mod {p^3}&\t{if $\sls{-1}p=\sls mp=-1$ and so $2p=x^2+my^2$.}\end{cases}\end{align*}}

\par\textbf{Theorem 1.12} {\sl Let $m\in\{3,5,11,29\}$ and
$F(m)=48^2,12^4,1584^2,396^4$ according as $m=3,5,11,29$. Suppose that $p$ is an odd prime such that $\sls{-2m}p=1$. Then

\begin{align*}\sum_{k=0}^{p-1}\f{\b{2k}k^2\b{4k}{2k}}{F(m)^k}
\e \begin{cases}4x^2-2p-\f{p^2}{4x^2}\mod {p^3}
\\\qq\t{if $\sls{-(-1)^{\f{m-1}2}2}p=\sls {(-1)^{\f{m-1}2}m}p=1$ and so $p=x^2+2my^2$,}
\\-8x^2+2p+\f{p^2}{8x^2}\mod {p^3}\\\qq\t{if $\sls{-(-1)^{\f{m-1}2}2}p=\sls {(-1)^{\f{m-1}2}m}p=-1$ and so $p=2x^2+my^2$.}\end{cases}\end{align*}}
\par Here the representability of $p$ and $2p$ by $x^2+my^2$, $x^2+2my^2$ or $2x^2+my^2$ is guaranteed by [24, Table 9.1].
\par It is noteworthy that the case of $p\equiv1\pmod{8}$ of Theorems 1.4-1.5 have been treated by Mao[11], and congruences in Theorems~1.6-1.7 in the weaker form under modulo~$p^{2}$ were proved by C. Wang and Z.W. Sun[30], and the congruence for $\sum_{k=0}^{p-1}\f{\b{2k}k^2\b{4k}{2k}}{81^k}$ modulo $p^2$ and the congruences modulo $p^2$ for the sums in Theorems 1.9-1.11 were first conjecturally formulated by Z.W. Sun[25,28]. In addition, the case $m=5$ in Theorem 1.11 has been treated by Beukers in [2, Corollary 1.19], while his result was stated in different form, and he did not give a proof of the corresponding rational CM value.

\par Combining Theorem 1.27 in Beukers' work [2] and the work [19] of the first author, in Section 5 we also confirm a number of conjectures of the first author given in [21] and [22] that can be stated as follows.
\v2
\par{\bf Theorem 1.13} {\it Let $p$ be a prime of the form $4k+1$ and so $p=x^2+4y^2$. Then
$$\Ls {-3}p\sum_{k=0}^{p-1}\f{\b{2k}k\b{3k}k\b{6k}{3k}}{12^{3k}}\e \Ls {33}p\sum_{ k=0}^{p-1}\f{\b{2k}k\b{3k}k\b{6k}{3k}}{66^{3k}}
\e 4x^2-2p-\f{p^2}{4x^2}\mod {p^3}.$$}

\par{\bf Theorem 1.14} {\it Let $p$ be a prime of the form $3k+1$ and so $p=x^2+3y^2$. Then
$$\sum_{k=0}^{p-1}\f{\b{2k}k\b{3k}k\b{6k}{3k}}{54000^k}
\e\Ls {5}p\Big(4x^2-2p-\f{p^2}{4x^2}\Big)\mod {p^3}.$$}

\par{\bf Theorem 1.15} {\it Let $p$ be a prime such that $p\e 1,3\mod 8$ and so $p=x^2+2y^2$. Then
$$\sum_{k=0}^{p-1}\f{\b{2k}k\b{3k}k\b{6k}{3k}}{20^{3k}}
\e\Ls {-5}p\Big(4x^2-2p-\f{p^2}{4x^2}\Big)\mod {p^3}.$$}

\par{\bf Theorem 1.16} {\it Let $p$ be an odd prime, $p\e 1,2,4\mod 7$ and so $p=x^2+7y^2$. Then
$$\Ls{-15}p\sum_{k=0}^{p-1}\f{\b{2k}k\b{3k}k\b{6k}{3k}}{(-15)^{3k}}
\e\Ls{-255}p\sum_{k=0}^{p-1}\f{\b{2k}k\b{3k}k\b{6k}{3k}}{255^{3k}}
\e 4x^2-2p-\f{p^2}{4x^2}\mod {p^3}.$$}

\par{\bf Theorem 1.17} {\it Let $p$ be a prime of the form $3k+1$ and so $4p=x^2+27y^2$. Then
$$\sum_{k=0}^{p-1}\f{\b{2k}k\b{3k}k\b{6k}{3k}}{(-12288000)^k}
\e\Ls {10}p\Big(x^2-2p-\f{p^2}{x^2}\Big)\mod {p^3}.$$}

\par{\bf Theorem 1.18} {\it Let $p$ be a prime such that $p\e 1,3,4,5,9\mod{11}$ and so $4p=x^2+11y^2$. Then
$$\sum_{k=0}^{p-1}\f{\b{2k}k\b{3k}k\b{6k}{3k}}{(-32)^{3k}}
\e\Ls {-2}p\Big(x^2-2p-\f{p^2}{x^2}\Big)\mod {p^3}.$$}

\par{\bf Theorem 1.19} {\it Let $p$ be a prime such that $\ls p{19}=1$ and so $4p=x^2+19y^2$. Then
$$\sum_{k=0}^{p-1}\f{\b{2k}k\b{3k}k\b{6k}{3k}}{(-96)^{3k}}
\e\Ls {-6}p\Big(x^2-2p-\f{p^2}{x^2}\Big)\mod {p^3}.$$}

\par{\bf Theorem 1.20} {\it Let $p$ be a prime such that $\ls p{43}=1$ and so $4p=x^2+43y^2$. Then
$$\sum_{k=0}^{p-1}\f{\b{2k}k\b{3k}k\b{6k}{3k}}{(-960)^{3k}}
\e\Ls {-15}p\Big(x^2-2p-\f{p^2}{x^2}\Big)\mod {p^3}.$$}

\par{\bf Theorem 1.21} {\it Let $p$ be a prime such that $\ls p{67}=1$ and so $4p=x^2+67y^2$. Then
$$\sum_{k=0}^{p-1}\f{\b{2k}k\b{3k}k\b{6k}{3k}}{(-5280)^{3k}}
\e\Ls {-330}p\Big(x^2-2p-\f{p^2}{x^2}\Big)\mod {p^3}.$$}

\par{\bf Theorem 1.22} {\it Let $p$ be a prime such that $\ls p{163}=1$ and so $4p=x^2+163y^2$. Then
$$\sum_{k=0}^{p-1}\f{\b{2k}k\b{3k}k\b{6k}{3k}}{(-640320)^{3k}}
\e\Ls {-10005}p\Big(x^2-2p-\f{p^2}{x^2}\Big)\mod {p^3}.$$}
\par It is noteworthy to remark that the congruences for the sums modulo $p^2$ in Theorems 1.18-1.22 were conjectured earlier by Z.W. Sun[28].

\par Next, we consider the Ap\'ery-like sequence  $\{V_n\}$ defined by
\begin{align*} V_n&=\sum_{k=0}^n\b
nk\b{n+k}k(-1)^k\b{2k}k^216^{n-k}=\sum_{k=0}^n\b{2k}k^2\b{2n-2k}{n-k}^2
\\& =\sum_{k=0}^n\b{2k}k^3\b k{n-k}(-16)^{n-k}.\end{align*}
In [23, Conjecture 22.29], the first author conjectured that for any prime $p>3$,
$$\sum_{n=0}^{p-1}\f{V_n}{8^n}
\e \sum_{n=0}^{p-1}\f{V_n}{(-16)^n} \e\begin{cases}
4x^2-2p-\f{p^2}{4x^2}\mod {p^3}&\t{if $p=x^2+4y^2\e 1\mod 4$,}
\\\f 34p^2\b{(p-3)/2}{(p-3)/4}^{-2}\mod {p^3}
&\t{if $p\e 3\mod 4$}\end{cases}$$
In Section 6, we shall verify the case of $p\equiv1\pmod{4}$ and obtain:
\v2
\par{\bf Theorem 1.23} {\it Let $p$ be a prime of the form $4k+1$ and so $p=x^2+4y^2$. Then
$$\sum_{n=0}^{p-1}\f{V_n}{8^n}
\e \sum_{n=0}^{p-1}\f{V_n}{(-16)^n} \e 4x^2-2p-\f{p^2}{4x^2}\mod {p^3}.$$}
\par We remark that  Appendix C in Beukers' paper [2] does not include the Ap\'ery-like sequence  $\{V_n\}$.
\par Let $\{T_n\}$ be the Ap\'ery-like sequence given by
$$T_n=\sum_{k=0}^n\b nk^2\b{2k}n^2=\sum_{k=0}^{[n/2]}\b{2k}k^2\b{4k}{2k}\b{n+2k}{4k}4^{n-2k}.$$
For an odd prime $p$, in [20] the first author proved that
$$ \sum_{n=0}^{p-1}\f{T_n}{4^n}
\e \begin{cases}
4x^2-2p-\f{p^2}{4x^2}\pmod {p^3} & \t{if $p=4k+1=x^2+4y^2$,}
\\-\f {p^2}4\b{\f{p-3}2}{\f{p-3}4}^{-2}\pmod
{p^3}&\t{if $p\e 3\mod 4$}\end{cases}$$
and conjectured the congruences for $\sum_{n=0}^{p-1}\f{T_n}{m^n}$ modulo $p^3$ with $m=1,-4,16$. In Section 7, using Beukers' theorem we confirm these conjectures stated as below.
\v2
\par{\bf Theorem 1.24} {\it Let $p$ be a prime such that
 $p\e 1,3\pmod 8$ and
so $p=x^2+2y^2$. Then
$$\sum_{n=0}^{p-1}\f{T_n}{(-4)^n}\e 4x^2-2p-\f {p^2}{4x^2}\pmod{p^3}.$$}
\par {\bf Theorem 1.25} {\it Let $p$ be an odd prime such that $p\e 1,2,4\pmod 7$ and so $p=x^2+7y^2$. Then
$$ \sum_{n=0}^{p-1}T_n
\e\sum_{n=0}^{p-1}\f{T_n}{16^n} \e
4x^2-2p-\f{p^2}{4x^2}\pmod{p^3}.$$}
\par In [20], the first author proved that for any prime $p\not=2,7$,
$$\sum_{n=0}^{p-1}T_n \e\begin{cases} 4x^2-2p\pmod{p^2}&\t{if $p\e 1,2,4\pmod
7$ and so $p=x^2+7y^2$,}
\\0\pmod {p^2}&\t{if $p\e 3,5,6\pmod 7$.}
\end{cases}$$
\par Let $\{D_n\}$ be the Domb numbers defined by
$$D_n=\sum_{k=0}^n\b nk^2\b{2k}k\b{2n-2k}{n-k} \q(n=0,1,2,\ldots).$$
In [21] the first author conjectured the congruences for $\sum_{n=0}^{p-1}\f{D_n}{m^n}$ modulo $p^3$ with $m\in\{1,-2,4,8,-8,16,-32,64\}$, where $p$ is an odd prime. The corresponding congruences modulo $p^2$ were conjectured by Z.W. Sun[28] earlier.
\par
In Section 8 we obtain the following results.
\v2
\par{\bf Theorem 1.26} {\it Let $p$ be a prime such that
 $p\e 1\pmod 8$ and
so $p=x^2+2y^2$. Then
$$\sum_{n=0}^{p-1}\f{D_n}{8^n}\e 4x^2-2p-\f {p^2}{4x^2}\pmod{p^3}.$$}
\par{\bf Theorem 1.27} {\it Let $p$ be a prime such that
 $p\e 1\pmod {12}$ and
so $p=x^2+3y^2$. Then
$$\sum_{n=0}^{p-1}\f{D_n}{(-2)^n}\e 4x^2-2p-\f {p^2}{4x^2}\pmod{p^3}.$$
Moreover, for $p\e 1\mod {24}$,
$$ \sum_{n=0}^{p-1}\f{D_n}{(-32)^n}\e 4x^2-2p-\f {p^2}{4x^2}\pmod{p^3}.$$}
\par{\bf Theorem 1.28} {\it Let $p$ be a prime such that $p\e 1,5\mod {24}$. Then
$$\sum_{n=0}^{p-1}\f{D_n}{(-8)^n}
\e\begin{cases} 4x^2-2p-\f{p^2}{4x^2}\mod{p^3}
\\\qq\qq\t{if $p\e 1\mod {24}$ and
so $p=x^2+6y^2$,}
\\8x^2-2p-\f{p^2}{8x^2}\mod {p^3}
\\\qq\qq\t{if $p\e 5\mod {24}$ and so $p=2x^2+3y^2$.}
\end{cases}$$}

\par Let $\{A_n\}$ be the Ap\'ery numbers $\{A_n\}$ given by
 $$ A_n= \sum_{k=0}^n\binom nk^2\binom{n+k}k^2\q(n=0,1,2,\ldots).$$
In Section 9, we shall prove the following congruence conjectured by the first author in [21].
\v2
 \par{\bf Theorem 1.29} {\it Let $p$ be a prime such that
 $p\e 1\pmod 3$ and
so $p=x^2+3y^2$. Then
$$\sum_{n=0}^{p-1}(-1)^nA_n\e 4x^2-2p-\f {p^2}{4x^2}\pmod{p^3}.$$}
\par In [27], Z.W. Sun conjectured the congruence for $\sum_{n=0}^{p-1}(-1)^nA_n$ modulo $p^2$ for any prime $p>3$.

\section*{2. Beukers' theorem and Weber's functions}

\par\q Let $\Bbb R$ be the set of real numbers and
$H=\{a+bi\mid a,b\in \Bbb R,\ b>0\}$. Let
$$
{\rm SL}_{2}(\mathbb{Z})=\left\{\left.\begin{pmatrix}
    a&b\\c&d
\end{pmatrix}\right|\, a,b,c,d\in\mathbb{Z},\,ad-bc=1\right\},
$$
and for any positive integer $N$,
$$
\Gamma_{0}(N)=\left\{\begin{pmatrix}
    a&b\\c&d
\end{pmatrix}\in{\rm SL}_{2}(\mathbb{Z})\mid c\equiv0\pmod{N}\right\}.
$$
For $\tau\in H$ the Dedekind eta function $\eta(\tau)$ is defined by
  $$\eta(\tau)=e^{2\pi i \tau/24}\prod_{n=1}^{\infty}\big(1-e^{2 \pi i\tau n}\big).$$
It is well known (see [7]) that
$$\eta(\tau+1)=e^{2\pi i/24}\eta(\tau)\q\t{and}\q \eta\big(-\f 1{\tau}\big)=\sqrt{-i\tau}\,\eta(\tau).$$
Furthermore, from [12] we know that
for any $\begin{pmatrix}
    a&b\\c&d
\end{pmatrix}\in{\rm SL}_{2}(\mathbb{Z})$,
$$ {\eta\left(\frac{a\tau+b}{c\tau+d}\right)} = \Big(\frac{a}{c_0}\Big)\zeta_{24}^{ab+cd(1-a^{2})-ca+3c_{0}(a-1)+r\frac{3}{2}(a^{2}-1)}
(c\tau+d)^{\frac{1}{2}}\eta(\tau),\eqno{(2.1)}$$
where
$\zeta_n=e^{2\pi i/n}$ and $c_0$ is given by $c=2^{r}c_{0}$ with $2\nmid c_0$.

\par For $c,d\in\Bbb Z$ let $(c,d)$ be the greatest common divisor of $c$ and $d$. For given imaginary
quadratic irrational number $\alpha$ let $\bar {\alpha}$ be the conjugate of $\alpha$. For an odd prime $p$ let $\Bbb Z_p$ be the ring of $p$-adic integers. Let $X_{0}(N)^{+}$ be the modular curve associated to the Fricke group $\Gamma_{0}(N)^{+}$, the Fuchsian group generated by $\Gamma_{0}(N)$ and its involution $\begin{pmatrix}  0&1\\-N&0 \end{pmatrix}$. In particular, recall that for a Fuchsian group $\Gamma$ commensurable with ${\rm SL}_{2}(\mathbb{Z})$ of genus zero, a Hauptmodul for $\Gamma$ is a modular function for $\Gamma$ with a unique pole.
\par In [2], Beukers established the following great theorem.
\v2
\par{\bf Beukers' theorem ([2, Theorems 1.15-1.16])} {\it Let $F(t)=\sum_{n=0}^{\infty}a_nt^n$.
 Suppose that
\par $\t{\rm (1)}$ $t(\tau)$ is a Hauptmodul for the modular group $\Gamma_0(N)^+$.
\par $\t{\rm (2)}$ $F(t(\tau))$ is a modular form of weight $2$ with respect to $\Gamma_0(N)^+$.
\par $\t{\rm (3)}$ $F(t(\tau))$ has a unique zero $($modulo the action of $\Gamma_0(N)^+)$ which is located at the pole of
$t(\tau)$. We assume it has order $1/r$ for some integer $r\ge 1$.
\par $\t{\rm (4)}$ $F(t(\tau))$ can be written as an $\eta$-product or $r\ge 4$.
\par Let $\alpha$ be an imaginary
quadratic number with positive imaginary part. Let $p$ be a prime not dividing $N$ which splits in the quadratic field
$\Bbb Q(\alpha)$.
 \par $\t{\rm (i)}$ Suppose that there exist integers $c$ and $d$ such that
\begin{align*}&(c,d)=1,\ N\mid c, \q c\alpha\bar{\alpha},c(\alpha+\bar{\alpha})\in\Bbb Z,
\\&p=(c\alpha+d)(c\bar{\alpha}+d)\q\t{and} \q c\alpha+d\notin p\Bbb Z_p.
\end{align*}
If $t(\alpha)\in\Bbb Z_p$, then
$$\sum_{n=0}^{p-1}a_nt(\alpha)^n\e (c\alpha+d)^2\mod {p^3}.$$
\par $\t{\rm (ii)}$ Suppose that there exist integers $c$ and $d$ such that
\begin{align*}&(c,dN)=1,\ Nc\alpha\bar{\alpha}\in\Bbb Z, \q c(\alpha+\bar{\alpha})\in\Bbb Z,
\\&p=N(c\alpha+d)(c\bar{\alpha}+d)\q\t{and} \q c\alpha+d\notin p\Bbb Z_p.
\end{align*}
If $t(\alpha)\in\Bbb Z_p$, then
$$\sum_{n=0}^{p-1}a_nt(\alpha)^n\e -N(c\alpha+d)^2\mod {p^3}.$$}

\par By checking Beukers' proof, Beukers' theorem(i) is also true when
$\Gamma_0(N)^+$ is replaced by $\Gamma_0(N)$.
\v2
\par For $\tau\in H$ let
$$g_2(\tau)=60\sum_{m,n\in\Bbb Z,(m,n)\not=(0,0)}\f 1{(m+n\tau)^4}
\q\t{and}\q\Delta(\tau)=(2\pi)^{12}\eta(\tau)^{24}.$$
The modular function $j(\tau)$ and its cubic root $\gamma_2(\tau)$ are defined by
$$j(\tau)=1728\f{g_2(\tau)^3}{\Delta(\tau)}\q\t{and}\q \gamma_2(\tau)
=12\f{g_2(\tau)}{\root 3\of {\Delta(\tau)}} .$$
Then $j(\tau)=\gamma_2(\tau)^3$.
Set $\omega=\f{-1+\sqrt{-3}}2$. It is well known (see [7, (12.4) and (12.6)]) that
\begin{align*}&\gamma_2(\tau+1)=\omega^2\gamma_2(\tau),\q \gamma_2\big(-\f 1{\tau}\big)=\gamma_2(\tau),
\\& \gamma_2\Big(\f{a\tau+b}{c\tau+d}\Big)=\omega^{ac-ab+a^2cd-cd}\gamma_2(\tau)
\q \t{for}\ a,b,c,d\in\Bbb Z\ \t{with}\ ad-bc=1\end{align*}
and so
\begin{align*}&j(\tau+1)=j(\tau),\q j\big(-\f 1{\tau}\big)=j(\tau),
\\& j\Big(\f{a\tau+b}{c\tau+d}\Big)=j(\tau)
\q \t{for}\ a,b,c,d\in\Bbb Z\ \t{with}\ ad-bc=1.\end{align*}
 For $\tau\in H$ let $f(\tau),f_1(\tau)$ and $f_2(\tau)$ be the Weber functions given by
\begin{align*}&f(\tau)=\zeta_{48}^{-1}\f{\eta((\tau+1)/2)}{\eta(\tau)}
=e^{-2\pi i\tau/48}\prod_{n=1}^{\infty}\big(1+e^{2\pi i\tau(n-\f 12)}\big),
\\&f_1(\tau)=\f{\eta(\tau/2)}{\eta(\tau)}=e^{-2\pi i\tau/48}
\prod_{n=1}^{\infty}\big(1-e^{2\pi i\tau(n-\f 12)}\big),
\\&f_2(\tau)=\sqrt 2
\f{\eta(2\tau)}{\eta(\tau)}=\sqrt 2\,e^{2\pi i\tau/24}\prod_{n=1}^{\infty}\big(1+e^{2\pi i\tau n}\big).\end{align*}
 From [7, (12.16), Corollary 12.19, Theorem 12.17] we know that
\begin{align*}&f(\tau)f_1(\tau)f_2(\tau)=f_1(2\tau)f_2(\tau)=\sqrt 2,
\\&f\Big(-\f 1{\tau}\Big)=f(\tau),\q f_1\Big(-\f 1\tau\Big)=f_2(\tau),\q f_2\Big(-\f 1\tau\Big)=f_1(\tau),
\\&\gamma_2(\tau)=\f{f(\tau)^{24}-16}{f(\tau)^8}=\f{f_1(\tau)^{24}+16}{f_1(\tau)^8}
=\f{f_2(\tau)^{24}+16}{f_2(\tau)^8},
\\&f(\tau+1)=\zeta_{48}^{-1}f_1(\tau),\q f_1(\tau+1)=\zeta_{48}^{-1}f(\tau),\q f_2(\tau+1)=\zeta_{24}f_2(\tau).\end{align*}
\par In order to prove our main results, we need the following basic lemmas.
\v2
\par{\bf Lemma 2.1} {\it
For two real numbers $a$ and $b$ with $b>0$ we have
\begin{align*}&\overline{f(a+\sqrt{-b})}=f(-a+\sqrt{-b}),\q \overline{f_1(a+\sqrt{-b})}=f_1(-a+\sqrt{-b}),\q \\&\overline{f_2(a+\sqrt{-b})}=f_2(-a+\sqrt{-b}).\end{align*}}
\par{\it Proof.} Note that
$$\overline{e^{2\pi i(a+\sqrt{-b})}}=\overline{e^{2\pi ia-2\pi \sqrt b}}=e^{-2\pi ia-2\pi \sqrt b}=e^{2\pi i(-a+\sqrt{-b})}.$$
We have
\begin{align*}\overline{f(a+\sqrt{-b})}
&=\overline{e^{-2\pi i(a+\sqrt{-b})/48}\prod_{n=1}^{\infty}\big(1+e^{2\pi i(a+\sqrt{-b})(n-\f 12)}\big)}
\\&=e^{-2\pi i(-a+\sqrt{-b})/48}\prod_{n=1}^{\infty}\big(1+e^{2\pi i(-a+\sqrt{-b})(n-\f 12)}\big)=f(-a+\sqrt{-b}),\end{align*}
\begin{align*}\overline{f_1(a+\sqrt{-b})}
&=\overline{e^{-2\pi i(a+\sqrt{-b})/48}\prod_{n=1}^{\infty}\big(1-e^{2\pi i(a+\sqrt{-b})(n-\f 12)}\big)}
\\&=e^{-2\pi i(-a+\sqrt{-b})/48}\prod_{n=1}^{\infty}\big(1-e^{2\pi i(-a+\sqrt{-b})(n-\f 12)}\big)=f_1(-a+\sqrt{-b}),
\\\overline{f_2(a+\sqrt{-b})}&=\overline{\sqrt 2\,e^{2\pi i(a+\sqrt{-b})/24}\prod_{n=1}^{\infty}\big(1+e^{2\pi i(a+\sqrt{-b}) n}\big)}
\\&=\sqrt 2\,e^{2\pi i(-a+\sqrt{-b})/24}\prod_{n=1}^{\infty}\big(1+e^{2\pi i(-a+\sqrt{-b}) n}\big)=f_2(-a+\sqrt{-b}).\end{align*}
This proves the lemma.
\v2
\par{\bf Lemma 2.2}
 {\it Let $\tau\in H$. Then the three roots of $x^3-\gamma_2(\tau)x+16=0$ are given by
 $x_0=-f(\tau)^8,\ x_1=f_1(\tau)^8$ and $x_2=f_2(\tau)^8$.  If $\gamma_2(\tau)^3\not=12^3$, then the three roots are distinct.}

\v2 \par{\it Proof}. From [31],
 $$f(\tau)^8=f_1(\tau)^8+f_2(\tau)^8.$$
Thus, $x_0+x_1+x_2=0$. Since $f(\tau)f_1(\tau)f_2(\tau)=\sqrt 2$ we have
$$x_0x_1x_2=-f(\tau)^8f_1(\tau)^8f_2(\tau)^8=-16.$$ On the other hand,
$$x_0x_1+x_0x_2+x_1x_2=x_0(-x_0)+\f{-16}{x_0}
=-\f{f(\tau)^{24}-16}{f(\tau)^8}=-\gamma_2(\tau).$$
Hence,
\begin{align*}&(x-x_0)(x-x_1)(x-x_2)\\&=x^3-(x_0+x_1+x_2)x^2
+(x_0x_1+x_0x_2+x_1x_2)x-x_0x_1x_2
\\&=x^3-\gamma_2(\tau)x+16.\end{align*}
This shows that $x_0,x_1$ and $x_2$ are the three roots of $x^3-\gamma_2(\tau)x+16=0$. The discriminant of $x^3-\gamma_2(\tau)x+16$
is given by $$D=-4(-\gamma_2(\tau))^3-27\cdot 16^2=4(\gamma_2(\tau)^3-12^3).$$
Hence, if $\gamma_2(\tau)^3\not=12^3$, then $D\not=0$ and so the three roots are pairwise distinct.

\v2
\par{\bf Lemma 2.3} {\it Let $p$ be an odd prime, and let $d>1$ be an integer not divisible by $p$. Suppose that $c$ is a positive integer, $cp=x^2+dy^2(x,y\in\Bbb Z)$ and $x+y\sqrt{-d}\not\in p\Bbb Z_p$. Then
 \begin{align*}&x+y\sqrt d\e 2x-\f {cp}{2x}-\f{c^2p^2}{8x^3}-\f{c^3p^3}{16x^5}\mod {p^4},
 \\&(x+y\sqrt{-d})^2\e 4x^2-2cp-\f{c^2p^2}{4x^2}-\f{c^3p^3}{8x^4}\mod {p^4}.
 \end{align*}}
  \par {\it Proof}. Suppose $A=x+y\sqrt {-d}$. Since $(A-x)^2=-dy^2=x^2-cp$, we have $A(A-2x)=-cp$
 and so $A\e 2x\mod p$. Set $A=2x+\f{kp}{2x}$. Then $(2x+\f{kp}{2x})\f{kp}{2x}=-cp$ and so
 $k+\f{k^2p}{4x^2}=-c.$  Hence
 \begin{align*} A&=2x+\f{kp}{2x}=2x-\f{(c+\f{k^2p}{4x^2})p}{2x}=2x-\f {cp}{2x}-\f{k^2p^2}{8x^3}
 \\&=2x-\f {cp}{2x}-\f{p^2}{8x^3}\Big(-c-\f{k^2p}{4x^2}\Big)^2
 \\&\e 2x-\f {cp}{2x}-\f{p^2}{8x^3}\Big(c^2+\f{ck^2p}{2x^2}\Big)
 \\&\e 2x-\f {cp}{2x}-\f{c^2p^2}{8x^3}-\f{c^3p^3}{16x^5}\mod {p^4}.\end{align*}
 Therefore,
  \begin{align*}A^2&\e \Big(2x-\f {cp}{2x}-\f{c^2p^2}{8x^3}-\f{c^3p^3}{16x^5}\Big)^2
 \\&\e \Big(2x-\f {cp}{2x}\Big)^2-2\Big(2x-\f {cp}{2x}\Big)\Big(\f{c^2p^2}{8x^3}+\f{c^3p^3}{16x^5}\Big)\end{align*}
 \begin{align*}&\e 4x^2-2cp+\f{c^2p^2}{4x^2}-\f{c^2p^2}{2x^2}+\f{c^3p^3}{8x^4}-\f{c^3p^3}{4x^4}
 \\&=4x^2-2cp-\f{c^2p^2}{4x^2}-\f{c^3p^3}{8x^4}\mod {p^4}.\end{align*}
 This completes the proof.
\v2
\section*{3. Proofs of congruences involving $\b{2k}k^3$}

\par{\bf Lemma 3.1} {\it For $\tau\in H$ let $t(\tau)$ be defined by
           $$ t(\tau)=\f q{\prod_{n=1}^{\infty}(1+q^{2n-1})^{24}}=\f 1{f(2\tau)^{24}},
    \q\t{where}\q
   q=e^{2\pi i\tau}.$$
    Then
    \begin{align*}&t(\tau)=\Ls{\eta(\tau)\eta(4\tau)}{\eta(2\tau)^2}^{24}
    =t\Big(-\f 1{4\tau}\Big),
    \\&(16t(\tau)-1)^3+j(2\tau)t(\tau)^2=0
    \end{align*}
    and $t(\tau)$ is a Hauptmodul for $\Gamma_{0}(4)^{+}$ with a unique pole at the cusp $[\frac{1}{2}]$.}
    \v2
\par{\it Proof.} It is clear that
\begin{align*}t(\tau)&=\f 1{f(2\tau)^{24}}
=\f{f_1(2\tau)^{24}f_2(2\tau)^{24}}{2^{12}}
=\f{\eta(\tau)^{24}}{\eta(2\tau)^{24}}
\cdot \f{\eta(4\tau)^{24}}{\eta(2\tau)^{24}}
\end{align*}
and
$$t(\tau)=\f 1{f(2\tau)^{24}}=\f 1{f(-\f 1{2\tau})^{24}}=t\Big(-\f 1{4\tau}\Big).$$
Since
$$j(2\tau)=\gamma_2(2\tau)^3=\f{(f(2\tau)^{24}-16)^3}{f(2\tau)^{24}},$$
we derive that

$$j(2\tau)=\f{(\f 1{t(\tau)}-16)^3}{\f 1{t(\tau)}}=\f{(1-16t(\tau))^3}{t(\tau)^2}$$
and so $(16t(\tau)-1)^3+j(2\tau)t(\tau)^2=0$.
\par
  Next we show that $t(\tau)$ is invariant under $\Gamma_{0}(4)$. For any $\begin{pmatrix}
     a&b\\c&d
 \end{pmatrix}\in \Gamma_{0}(4)$, one can have that
 \begin{align*}
     t\Big(\frac{a\tau+b}{c\tau+d}\Big)
     =\Ls{\eta\big(\frac{a\tau+b}{c\tau+d}\big)
     \eta\big(4\frac{a\tau+b}{c\tau+d}\big)}
     {\eta\big(2\frac{a\tau+b}{c\tau+d}\big)^2}^{24}
     =\Ls{\eta\big(\frac{a\tau+b}{c\tau+d}\big)\eta
     \big(\frac{a(4\tau)+4b}{(c/4)(4\tau)+d}\big)}
     {\eta\big(\frac{a(2\tau)+2b}{(c/2)(2\tau)+d}\big)^2}^{24}.
 \end{align*}
 By (2.1), one can find that
 \begin{align*}
     \eta\Big(\frac{a\tau+b}{c\tau+d}\Big)^{24}&=(c\tau+d)^{12}\eta(\tau)^{24},\\
     \eta\Big(\frac{a(4\tau)+4b}{(c/4)(4\tau)+d}\Big)^{24}&=((c/4)(4\tau)+d)^{12}\eta(4\tau)=(c\tau+d)^{12}\eta(4\tau)^{24},\\
     \eta\Big(\frac{a(2\tau)+2b}{(c/2)(2\tau)+d}\Big)^{48}&=((c/2)(2\tau)+d)^{24}\eta(2\tau)=(c\tau+d)^{24}\eta(2\tau)^{48}
 \end{align*}
 and thus deduce that
 $$
 t\Big(\frac{a\tau+b}{c\tau+d}\Big)=t(\tau)
 \q\t{for any $\begin{pmatrix}
     a&b\\c&d
 \end{pmatrix}\in \Gamma_{0}(4)$}.$$ Combining this with the relation $t(\tau)=t(-\f 1{4\tau})$ proven above, one can tell that $t(\tau)$ is invariant under $\Gamma_{0}(4)^{+}$. To prove that $t(\tau)$ is a Hauptmodul for $\Gamma_{0}(4)^{+}$, it suffices to show that $t(\tau)$ has a unique pole in the modular curve $X_{0}(4)^{+}$. To this end, one first notices that $t(\tau)$ has no zeros or poles in $H$ by its infinite product representation, and $t(\tau)=q+O(q^{2})$ has a simple zero at the cusp $[i\infty]$ by definition. Since $X_{0}(4)^{+}$ has only two cusps, $[i\infty]$ and $[\frac{1}{2}]$, and $t(\tau)$ is a modular function on $X_{0}(4)^{+}$ with no zeros or poles in $H$, then by the Residue Theorem for compact Riemann surfaces, one can see that this forces $t(\tau)$ to have a simple pole at the cusp $[\frac{1}{2}]$ which is the unique pole of $t(\tau)$ in $X_{0}(4)^{+}$. Hence, $t(\tau)$ is a Hauptmodul for $\Gamma_{0}(4)^{+}$.
This completes the proof.
\v2

 \par{\bf Lemma 3.2} {\it For $\tau\in H$ let $t(\tau)$ be given in Lemma 3.1. Then
    \begin{align*}&
     t\Ls{3+\sqrt{-7}}{8}=1,\ t\Ls{\sqrt{-7}}2=\f 1{4096},
    \\&t\Ls{3+\sqrt{-3}}4=\f 1{16},\ t\Ls{\sqrt{-3}}2=\f 1{256},
    \\&t\Ls{1+i}2=-\f 18,\q t\Ls{1+\sqrt{-2}}2=-\f 1{64}.
    \end{align*}}
    \par{\it Proof.}  From [7, (12.20)] we have  $j(\f{3+\sqrt{-7}}2)=-15^3$.
  Since $j(\tau+1)=j(\tau)$ and $j(\tau)=j(-\f 1{\tau})$, we see that
  \begin{align*}j\Ls{3+\sqrt{-7}}4&=j\Big( \f{-1+\sqrt{-7}}4\Big)
  =j\Big(-\f 4{-1+\sqrt{-7}}\Big)\\&=j\Ls{1+\sqrt{-7}}2
  =j\Ls{3+\sqrt{-7}}2=-15^3.\end{align*}
   Hence, appealing to Lemma 3.1 we see that $t(\f{3+\sqrt{-7}}8)$ is a root of $(16x-1)^3-15^3x^2=(x-1)(4096x^2-47x+1)=0$ and so
   $$t\Big(\f{3+\sqrt{-7}}8\Big)\in\Big\{1,\f{47+45\sqrt{-7}}{8192},
   \f{47-45\sqrt{-7}}{8192}\Big\}.$$
  By Lemma 2.1 and $t(\tau)=t(-\f 1{4\tau})$ given in Lemma 3.1,
\begin{align*}
    \overline{t\Big(\frac{3+\sqrt{-7}}{8}\Big)}
    &=\overline{f\Ls{3+\sqrt{-7}}4^{-24}}=
    \overline{f\Ls{3+\sqrt{-7}}4}^{-24}=f\Ls{-3+\sqrt{-7}}4^{-24}
    \\&=t\Big(\frac{-3+\sqrt{-7}}{8}\Big)=t\Big(-\frac{1}{4\cdot
    \frac{-3+\sqrt{-7}}{8}}\Big)=t\Big(\frac{3+\sqrt{-7}}{8}\Big).
\end{align*}
Thus, $t\big(\frac{3+\sqrt{-7}}{8}\big)$ is real and therefore $t\sls{3+\sqrt{-7}}8=1$.
\par From [7, (12.20)], $\gamma_2(\sqrt{-7})=255$. By Lemma 2.2,
 $-f(\sqrt{-7})^8,f_1(\sqrt{-7})^8$ and $f_2(\sqrt{-7})^8$ are the distinct roots of $x^3-255x+16
=(x+16)((x-8)^2-63)=0$ and so
$$-f(\sqrt{-7})^8,f_1(\sqrt{-7})^8,f_2(\sqrt{-7})^8\in\big\{-16,8+3\sqrt 7,8-3\sqrt 7\big\}.$$
Hence
 $$t\Ls{\sqrt{-7}}2=\f 1{f(\sqrt{-7})^{24}}\in\Big\{\f 1{16^3},-\f 1{(8+3\sqrt 7)^3},
 -\f 1{(8-3\sqrt 7)^3}\Big\}.$$
 Since $e^{2\pi i\sqrt{-7}/2}=e^{-\sqrt 7\pi}$ and so
 $$t\Ls{\sqrt{-7}}2 =e^{-\sqrt 7\pi}\prod_{n=1}^{\infty}\Big(1+e^{-\sqrt 7\pi(2n-1)}\Big)^{-24}>0,$$
 we must have $t\ls{\sqrt{-7}}2=\f 1{16^3}=\f 1{4096}$.

     \par Since $j(\tau+1)=j(\tau)$, from [7, p.261] we have $j\ls{3+\sqrt{-3}}2=j\ls{1+\sqrt{-3}}2=0$ and so $t\ls{3+\sqrt{-3}}4=\f 1{16}$ by Lemma 3.1.
 \par From [7, pp.261,291] we have $j(\sqrt{-3})=54000$. Thus, applying Lemma 3.1 we see that $t(\f{\sqrt{-3}}2)$ is a root of $(16x-1)^3+54000x^2=(256x-1)(16x^2+208x+1)=0$ and so
 $$t\Ls{\sqrt{-3}}2\in\Big\{\f 1{256},\f{-26+15\sqrt 3}4,\f{-26-15\sqrt 3}4\Big\}.$$
 Since $e^{2\pi i\sqrt{-3}/2}=e^{-\sqrt 3\pi}$ and so
 $$t\Ls{\sqrt{-3}}2=e^{-\sqrt 3\pi}\prod_{n=1}^{\infty}\Big(1+e^{-\sqrt 3\pi(2n-1)}\Big)^{-24}>0,$$
 we must have $t\ls{\sqrt{-3}}2=\f 1{256}$.

    \par  By [7, (12.20)], $j(i)=12^3$ and so $j(1+i)=12^3$ since $j(\tau+1)=j(\tau)$. From Lemma 3.1, $t\sls{1+i}2$ is a root of $(16x-1)^3+12^3x^2=(8x+1)^2(64x-1)=0$. Thus, $t\sls{1+i}2=-\f 18$ or $\f 1{64}$. Note that $e^{2\pi i\f{1+i}2}=-e^{-\pi}$. We see that
 $$t\Ls{1+i}2=\f{-e^{-\pi}}{\prod_{n=1}^{\infty}
 (1-e^{-(2n-1)\pi})^{24}}<0$$
 and so $t\sls{1+i}2=-\f 18$.
  \par By [7, Corollary 12.19], $f(\tau+1)=e^{-\f{2\pi i}{48}}f_1(\tau)$. Thus, $f(\tau+1)^{24}=e^{-\pi i}f_1(\tau)^{24}=-f_1(\tau)^{24}$.
    Since $f_2(-\f 1{\tau})=f_1(\tau)$ and $f_1(2\tau)f_2(\tau)=\sqrt 2$, we see that
$f_2\sls{\sqrt{-2}}2=f_2(-\f 1{\sqrt{-2}})=f_1(\sqrt {-2})$
and so $f_1(\sqrt{-2})^2=f_1(\sqrt {-2})f_2\sls{\sqrt{-2}}2=\sqrt 2$.
  Therefore,
  $$t\Ls{1+\sqrt{-2}}2=\f 1{f(1+\sqrt{-2})^{24}}=-\f 1{f_1(\sqrt{-2})^{24}}=-\f 1{(\sqrt 2)^{12}}=-\f 1{64}.$$
    This completes the proof.
\v2
\par{\bf Remark 3.1} Lemma 3.2 was stated by Beukers in [2, p.29] without proof.
 \v2
\par{\bf Lemma 3.3 } {\it Let $t(\tau)$ be given in Lemma 3.1. For $\tau\in H$ we have
\begin{align*}
        \sum_{n=0}^{\infty}\binom{2n}{n}^{3}t(\tau)^{n}
        =\frac{\eta(2\tau)^{20}}{\eta(\tau)^{8}\eta(4\tau)^{8}}
           \end{align*}
is a weight $2$ modular form for $\Gamma_{0}(4)^{+}$ with a unique zero at the cusp $[\frac{1}{2}]$ of order $1$.}
\v2
\par{\it Proof.} By Lemma 3.1, $t(\tau)=\ls{\eta(\tau)\eta(4\tau)}{\eta(2\tau)^2}^{24}$. Thus, the equality in Lemma 3.3 was already proved by Cooper in [6, Theorem 4.1(d)]. It remains to justify the eta quotient on the right hand side to be a weight~2 modular form for $\Gamma_{0}(4)^{+}$ with a unique zero at the cusp $[\frac{1}{2}]$ of order~1. To this end, note by (2.1) that for any $\begin{pmatrix}
    a&b\\c&d
\end{pmatrix}\in\Gamma_{0}(4)$,
\begin{align*}
      \eta\Big(\frac{a\tau+b}{c\tau+d}\Big)^8&= \omega^{ab+cd(1-a^{2})-ca+3c_{0}(a-1)+\frac{3}{2}(a^{2}-1)r}
(c\tau+d)^{4}\eta(\tau)^{8},\\
\eta\Big(\frac{a(2\tau)+2b}{(c/2)(2\tau)+d}\Big)^{20}
&=\omega^{\frac{5}{2}(a(2b)+(c/2)d(1-a^{2})
-\frac{1}{2}ca+3c_{0}(a-1)+\frac{3}{2}(a^{2}-1))(r-1)}
\\&\q\times((c/2)(2\tau)+d)^{10}\eta(2\tau)^{20},\\
\eta\Big(\frac{a(4\tau)+4b}{(c/4)(4\tau)+d}\Big)^{8}
&=\omega^{a(4b)+(c/4)d(1-a^{2})-\frac{1}{4}ca+3c_{0}(a-1)+\frac{3}{2}(a^{2}-1))(r-2)}
\\&\q\times((c/4)(4\tau)+d)^{4}\eta(4\tau)^{8}.
\end{align*}
Then after some cancellation, one deduces that
$$
\frac{\eta(\frac{a(4\tau)+4b}{(c/4)(4\tau)+d})^{20}}{\eta(\frac{a\tau+b}{c\tau+d})^8\eta(\frac{a(4\tau)+4b}{(c/4)(4\tau)+d})^{8}}=(c\tau+d)^{2}\frac{\eta(2\tau)^{20}}{\eta(\tau)^{8}\eta(4\tau)^{8}}.
$$
Similarly, under the action of the Fricke involution $\begin{pmatrix}
    0&-1\\4&0
\end{pmatrix}$,
\begin{align*}
    \eta(-1/4\tau)^{8}&=(4\tau)^{4}\eta(4\tau)^{8},\\
    \eta(-1/2\tau)^{20}&=(2\tau)^{10}\eta(2\tau)^{20},\\
    \eta(-1/\tau)^{8}&=\tau^{4}\eta(\tau)^{8},
\end{align*}
and thus,
$$
\frac{\eta(-1/2\tau)^{20}}{\eta(-1/4\tau)^{8}\eta(-1/\tau)^{8}}=(2\tau)^{2}\frac{\eta(2\tau)^{20}}{\eta(\tau)^{8}\eta(4\tau)^{8}}.
$$
Therefore, $\frac{\eta(2\tau)^{20}}{\eta(\tau)^{8}\eta(4\tau)^{8}}$ is a weight~2 modular form for $\Gamma_{0}(4)^{+}$. Moreover, by its infinite product representation, it is clear that $\frac{\eta(2\tau)^{20}}{\eta(\tau)^{8}\eta(4\tau)^{8}}$ has no zeros or poles in $H$, and $\frac{\eta(2\tau)^{20}}{\eta(\tau)^{8}\eta(4\tau)^{8}}=1+O(q)$ is of order~0 at the cusp $[i\infty]$. So since $X_{0}(4)^{+}$ has only two cusps $[i\infty]$ and $[\frac{1}{2}]$, and  $\frac{\eta(2\tau)^{20}}{\eta(\tau)^{8}\eta(4\tau)^{8}}$ is a weight~2 modular form for $\Gamma_{0}(4)^{+}$ with no zeros or poles apart from the cusp $[\frac{1}{2}]$, then by the Riemann-Roch theorem this forces $\frac{\eta(2\tau)^{20}}{\eta(\tau)^{8}\eta(4\tau)^{8}}$ to have a zero of order~1 at the cusp $[\frac{1}{2}]$. This completes the proof.

\v2
\par Let $t(\tau)$ be given in Lemma 3.1.
From Lemmas 3.1 and 3.3 one can see that both
$t(\tau)$ and $\sum_{n=0}^{\infty}\binom{2n}{n}^{3}t(\tau)^{n}
$ satisfy the assumptions in Beukers' theorem.\v2
\par{\bf Proof of Theorem 1.1}. Since $p=x^2+7y^2=(x+y\sqrt{-7})(x-y\sqrt{-7})$ for $x,y\in\Bbb Z$, one may choose the sign of $y$ so that $x+y\sqrt{-7}\notin p\Bbb Z_p$.
 Set
 $$a_n=\b{2n}n^3,\q c=8y,\q d=x-3y,\q \alpha=\f{3+\sqrt{-7}}8\q\t{and}\q N=4.$$
 Then clearly
$(c,d)=1,\ N\mid c, \q c\alpha\bar{\alpha},c(\alpha+\bar{\alpha})\in\Bbb Z,$ $p=(c\alpha+d)(c\bar{\alpha}+d)$ and $c\alpha+d\notin p\Bbb Z_p$.
 By Lemma 3.2, $t\sls{3+\sqrt{-7}}8=1$.
 Hence, applying Beukers' theorem(i) and Lemma 2.3 we obtain
 \begin{align*}\sum_{n=0}^{p-1}\b{2n}n^3
 &=\sum_{n=0}^{p-1}\b{2n}n^3t\Ls{3+\sqrt{-7}}8^n\e \Big(8y\cdot \f{3+\sqrt{-7}}8+x-3y\Big)^2
 \\&=(x+y\sqrt {-7})^2\e 4x^2-2p-\f{p^2}{4x^2}\mod {p^3}.\end{align*}
\par For $p=x^2+7y^2\e 1\mod 4$ we have $2\mid y$.
Set $$a_n=\b{2n}n^3,\q
c=2y,\q d=x,\q \alpha=\f{\sqrt{-7}}2\q\t{and}\q N=4.$$
Then $(c,d)=1,\ N\mid c, \q c\alpha\bar{\alpha},c(\alpha+\bar{\alpha})\in\Bbb Z,$ $p=(c\alpha+d)(c\bar{\alpha}+d)$ and $c\alpha+d\notin p\Bbb Z_p$.
Since $t\ls{\sqrt{-7}}2=\f 1{4096}$ by Lemma 3.2, using Beukers' theorem(i) and Lemma 2.3 we see that
 \begin{align*}\sum_{n=0}^{p-1}\f{\b{2n}n^3}{4096^n}
 &=\sum_{n=0}^{p-1}\b{2n}n^3t\Ls{\sqrt{-7}}2^n\e \Big(2y\cdot \f{\sqrt{-7}}2+x\Big)^2
 \\&\e 4x^2-2p-\f{p^2}{4x^2}\mod {p^3}.\end{align*}

\par For $p=x^2+7y^2\e 3\mod 4$ we have $2\mid x$ and $2\nmid y$.
Set $$a_n=\b{2n}n^3,\q
c=y,\q d=\f x2,\q \alpha=\f{\sqrt{-7}}2\q\t{and}\q N=4.$$
Then
$(c,dN)=1,\ Nc\alpha\bar{\alpha}\in\Bbb Z, \q c(\alpha+\bar{\alpha})\in\Bbb Z,\
p=N(c\alpha+d)(c\bar{\alpha}+d)$ and $c\alpha+d\notin p\Bbb Z_p$.
Since $t\ls{\sqrt{-7}}2=\f 1{4096}$, using Beukers' theorem(ii) and Lemma 2.3 we see that
 \begin{align*}\sum_{n=0}^{p-1}\f{\b{2n}n^3}{4096^n}
 &=\sum_{n=0}^{p-1}\b{2n}n^3t\Ls{\sqrt{-7}}2^n\e
 -4\Big(y\cdot\f{\sqrt{-7}}2+\f x2\Big)^2\\&=-(x+y\sqrt {-7})^2\e -\Big(4x^2-2p-\f{p^2}{4x^2}\Big)\mod {p^3}.\end{align*}
\par Note that $p\mid \b{2k}k$ and so $p^3\mid \b{2k}k^3$ for $k=\f{p+1}2,\ldots,p-1$.
Summarizing the above proves Theorem 1.1.
\v2

\par{\bf Proof of Theorem 1.2}. Since $p=x^2+3y^2=(x+y\sqrt{-3})(x-y\sqrt{-3})$ for $x,y\in\Bbb Z$, one may choose the sign of $y$ so that $x+y\sqrt{-3}\notin p\Bbb Z_p$.
 Set
 $$a_n=\b{2n}n^3,\q c=4y,\q d=x-3y,\q \alpha=\f{3+\sqrt{-3}}4\q\t{and}\q N=4.$$
 Then clearly
 $(c,d)=1,\ N\mid c,$ $c\alpha\bar{\alpha},\ c(\alpha+\bar{\alpha})\in\Bbb Z$,
$p=(c\alpha+d)(c\bar{\alpha}+d)$ and $c\alpha+d\notin p\Bbb Z_p$.
 By Lemma 3.2, $t\sls{3+\sqrt{-3}}4=\f 1{16}$.
 Hence, applying Beukers' theorem(i) and Lemma 2.3 we obtain
 \begin{align*}\sum_{n=0}^{p-1}\f{\b{2n}n^3}{16^n}
 &=\sum_{n=0}^{p-1}\b{2n}n^3t\Ls{3+\sqrt{-3}}4^n\e \Big(4y\cdot \f{3+\sqrt{-3}}4+x-3y\Big)^2
 \\&=(x+y\sqrt {-3})^2\e 4x^2-2p-\f{p^2}{4x^2}\mod {p^3}.\end{align*}
\par For $p=x^2+3y^2\e 1\mod 4$ we have $2\mid y$.
Set $$a_n=\b{2n}n^3,\q
c=2y,\q d=x,\q \alpha=\f{\sqrt{-3}}2\q\t{and}\q N=4.$$
Then $(c,d)=1,\ N\mid c, \q c\alpha\bar{\alpha},$ $c(\alpha+\bar{\alpha})\in\Bbb Z,$
$p=(c\alpha+d)(c\bar{\alpha}+d)$ and $c\alpha+d\notin p\Bbb Z_p$.
Since $t\ls{\sqrt{-3}}2=\f 1{256}$ by Lemma 3.2, using Beukers' theorem(i) and Lemma 2.3 we see that
 \begin{align*}\sum_{n=0}^{p-1}\f{\b{2n}n^3}{256^n}
 &=\sum_{n=0}^{p-1}\b{2n}n^3t\Ls{\sqrt{-3}}2^n\e \Big(2y\cdot \f{\sqrt{-3}}2+x\Big)^2
 \\&\e 4x^2-2p-\f{p^2}{4x^2}\mod {p^3}.\end{align*}
\par For $p=x^2+3y^2\e 3\mod 4$ we have $2\mid x$ and $2\nmid y$.
Set $$a_n=\b{2n}n^3,\q
c=y,\q d=\f x2,\q \alpha=\f{\sqrt{-3}}2\q\t{and}\q N=4.$$
Then $(c,dN)=1,\ Nc\alpha\bar{\alpha}\in\Bbb Z,$  $c(\alpha+\bar{\alpha})\in\Bbb Z$,
$p=N(c\alpha+d)(c\bar{\alpha}+d)$ and $c\alpha+d\notin p\Bbb Z_p$.
Since $t\ls{\sqrt{-3}}2=\f 1{256}$ by Lemma 3.2, using Beukers' theorem(ii) and Lemma 2.3 we see that
 \begin{align*}\sum_{n=0}^{p-1}\f{\b{2n}n^3}{256^n}
 &=\sum_{n=0}^{p-1}\b{2n}n^3 t\Ls{\sqrt{-3}}2^n\e
 -4\Big(y\cdot\f{\sqrt{-3}}2+\f x2\Big)^2\\&=-(x+y\sqrt {-3})^2\e -\Big(4x^2-2p-\f{p^2}{4x^2}\Big)\mod {p^3}.\end{align*}
 \par Note that $p\mid \b{2k}k$ and so $p^3\mid \b{2k}k^3$ for $k=\f{p+1}2,\ldots,p-1$. From the above we deduce Theorem 1.2.
\v2
\par{\bf Proof of Theorem 1.3}. Since $p=(x+2yi)(x-2yi)$, we may choose the sign of $x$ so that $x+2yi\not\in p\Bbb Z_p$.  Set
 $$a_n=\b{2n}n^3,\q c=4y,\q d=x-2y,\q \alpha=\f{1+i}2\q\t{and}\q N=4.$$
 Then clearly
 $(c,d)=1,\ N\mid c,$ $c\alpha\bar{\alpha},\ c(\alpha+\bar{\alpha})\in\Bbb Z$,
$p=(c\alpha+d)(c\bar{\alpha}+d)$ and $c\alpha+d\notin p\Bbb Z_p$.
 By Lemma 3.2, $t\sls{1+i}2=-\f 18$.
 Hence, applying Beukers' theorem(i) and Lemma 2.3 we obtain
 \begin{align*}\sum_{n=0}^{p-1}\f{\b{2n}n^3}{(-8)^n}
 &=\sum_{n=0}^{p-1}\b{2n}n^3t\Ls{1+i}2^n\e \Big(4y\cdot \f{1+i}2+x-2y\Big)^2
 \\&=(x+y\sqrt{-4})^2\e 4x^2-2p-\f{p^2}{4x^2}\mod {p^3}.\end{align*}
 Since $p\mid \b{2k}k$ and so $p^3\mid \b{2k}k^3$ for $k=\f{p+1}2,\ldots,p-1$, the result follows.
\v2
\par{\bf Proof of Theorem 1.4}. Since $p=x^2+2y^2=(x+y\sqrt{-2})(x-y\sqrt{-2})$ for $x,y\in\Bbb Z$, one may choose the sign of $y$ so that $x+y\sqrt{-2}\notin p\Bbb Z_p$.
\par For $p=x^2+2y^2\e 1\mod 8$ we have $2\mid y$.
Set $$a_n=\b{2n}n^3,\q
c=2y,\q d=x-y,\q \alpha=\f{1+\sqrt{-2}}2\q\t{and}\q N=4.$$
Then $(c,d)=1,\ N\mid c, \q c\alpha\bar{\alpha},$ $c(\alpha+\bar{\alpha})\in\Bbb Z,$
$p=(c\alpha+d)(c\bar{\alpha}+d)$ and $c\alpha+d\notin p\Bbb Z_p$.
Since $t\ls{1+\sqrt{-2}}2=-\f 1{64}$ by Lemma 3.2, using Beukers' theorem(i) and Lemma 2.3 we see that
 \begin{align*}\sum_{n=0}^{p-1}\f{\b{2n}n^3}{(-64)^n}
 &=\sum_{n=0}^{p-1}\b{2n}n^3t\Ls{1+\sqrt{-2}}2^n\e \Big(2y\cdot \f{1+\sqrt{-2}}2+x-y\Big)^2
 \\&=(x+y\sqrt {-2})^2\e 4x^2-2p-\f{p^2}{4x^2}\mod {p^3}.\end{align*}

\par For $p=x^2+2y^2\e 3\mod 8$ we have $2\nmid xy$.
Set $$a_n=\b{2n}n^3,\q
c=y,\q d=\f{x-y}2,\q \alpha=\f{1+\sqrt{-2}}2\q\t{and}\q N=4.$$
Then $(c,dN)=1,\ Nc\alpha\bar{\alpha}\in\Bbb Z,$  $c(\alpha+\bar{\alpha})\in\Bbb Z$,
$p=N(c\alpha+d)(c\bar{\alpha}+d)$ and $c\alpha+d\notin p\Bbb Z_p$.
Since $t\ls{1+\sqrt{-2}}2=-\f 1{64}$ by Lemma 3.2, using Beukers' theorem(ii) and Lemma 2.3 we see that
 \begin{align*}\sum_{n=0}^{p-1}\f{\b{2n}n^3}{(-64)^n}
 &=\sum_{n=0}^{p-1}\b{2n}n^3 t\Ls{1+\sqrt{-2}}2^n\e
 -4\Big(y\cdot\f{1+\sqrt{-2}}2+\f{x-y}2\Big)^2\\&=-(x+y\sqrt {-2})^2\e -\Big(4x^2-2p-\f{p^2}{4x^2}\Big)\mod {p^3}.\end{align*}
 \par Note that $p\mid \b{2k}k$ and so $p^3\mid \b{2k}k^3$ for $k=\f{p+1}2,\ldots,p-1$. From the above Theorem 1.4 is proved.

\section*{4. Proofs of congruences involving $\b{2k}k^2\b{4k}{2k}$}

\v2
 \par For given rational number $n$ set
 $$\sigma(n)=\begin{cases} \sum_{d\mid n}d&\t{if $n\in\{1,2,3,\ldots\}$,}
 \\0&\t{otherwise.}\end{cases}$$
 For $\tau\in H$ let $u(\tau)$ be defined by
$$
u(\tau)=\frac{f_{2}(\tau)^{24}}{(f_{2}(\tau)^{24}+64)^{2}}.
\eqno{(4.1)}
$$
Since $f_1(2\tau)f_2(\tau)=\sqrt 2$, we see that
$$u(\tau)=\f{f_2(\tau)^{24}}{(f_2(\tau)^{24}+64)^2}
=\f{2^{12}/f_1(2\tau)^{24}}{(2^{12}/f_1(2\tau)^{24}+64)^2}
=\f{f_1(2\tau)^{24}}{(f_1(2\tau)^{24}+64)^2}.\eqno{(4.2)}$$
As $f_2(\tau)=f_1(-\f 1{\tau})$, we also have
$$u(\tau)=\f{f_1(-\f 1{\tau})^{24}}{(f_1(-\f 1{\tau})^{24}+64)^2}=u\Big(-\f 1{2\tau}\Big).\eqno{(4.3)}$$
 \par{\bf Lemma 4.1} {\it For $\tau\in H$ let $u(\tau)$ be given in (4.1) and $E_{2}(\tau)$ denote the normalized weight $2$ Eisenstein series
$$
E_{2}(\tau)=1-24\sum_{n=1}^{\infty}\frac{nq^{n}}{1-q^{n}}=1-24\sum_{n=1}^{\infty}
\sigma(n)q^n.
$$
Then one has that
$$
u(\tau)=\left(\frac{\eta(\tau)^{2}\eta(2\tau)^{2}}{8E_{2}(2\tau)-4E_{2}(\tau)}\right)^{4}
$$
is a Hauptmodul for $\Gamma_{0}(2)^{+}$ with a unique pole at $\frac{1+i}{2}$,
and for $\tau\in H$ such that $u(\tau)$ is near~0,
$$
\sum_{n=0}^{\infty}\binom{2n}n^2\binom{4n}{2n}u(\tau)^n
=2E_{2}(2\tau)-E_{2}(\tau)=1+24\sum_{n=1}^{\infty}\Big(\sigma(n)-2\sigma
\ls n2\Big)q^n$$
is a weight $2$ modular form with respect to $\Gamma_{0}(2)^{+}$ with a unique zero of order~$\frac{1}{4}$ at $\frac{1+i}{2}$.
}
\v2
\par {\it Proof}. Both identities were proved by Cooper in [5, Theorem 4.30]. For any $\begin{pmatrix}
    a&b\\c&d
\end{pmatrix}\in\Gamma_{0}(2)$, write
\begin{align*}
    u\Big(\frac{a\tau+b}{c\tau+d}\Big)
    =\Big(\frac{\eta(\frac{a\tau+b}{c\tau+d})^{2}
    \eta(2\frac{a\tau+b}{c\tau+d})^{2}}{8E_{2}(2\frac{a\tau+b}
    {c\tau+d})-4E_{2}(\frac{a\tau+b}{c\tau+d})}\Big)^{4}
    =\Big(\frac{\eta(\frac{a\tau+b}{c\tau+d})^{2}
    \eta(\frac{a(2\tau)+2b}{(c/2)(2\tau)+d})^{2}}{8E_{2}
    (\frac{a(2\tau)+b}{(c/2)(2\tau)+d})-4E_{2}(\frac{a\tau+b}{c\tau+d})}
    \Big)^{4}.
\end{align*}
Aside of this, by the transformation formula (see [4, Corollary 5.2.17])
$$
E_{2}\left(\frac{a\tau+b}{c\tau+d}\right)=(c\tau+d)^{2}E_{2}(\tau)
-\frac{6ic}{\pi}(c\tau+d)
\q\t{for $\begin{pmatrix}
    a&b\\c&d
\end{pmatrix}\in{\rm SL}_{2}(\mathbb{Z})$},$$
 one can deduce that
\begin{align*}
    &8E_{2}\Big(\frac{a(2\tau)+b}{(c/2)(2\tau)+d}\Big)
    -4E_{2}\Big(\frac{a\tau+b}{c\tau+d}\Big)\\
    &=8(c\tau+d)^{2}E_{2}(2\tau)-\frac{24ic}{\pi}(c\tau+d)-\left(4(c\tau+d)^{2}E_{2}(\tau)-\frac{24ic}{\pi}(c\tau+d)\right)\\
    &=(c\tau+d)^{2}\left(8E_{2}(2\tau)-4E_{2}(\tau)\right).
\end{align*}
Also, by (2.1), one has that
\begin{align*}
      \eta\Big(\frac{a\tau+b}{c\tau+d}\Big)^8&= \omega^{ab+cd(1-a^{2})-ca+3c_{0}(a-1)+\frac{3}{2}(a^{2}-1)r}
(c\tau+d)^{4}\eta(\tau)^{8},\\
\eta\Big(\frac{a(2\tau)+2b}{(c/2)(2\tau)+d}\Big)^8
&=\omega^{a(2b)+(c/2)d(1-a^{2})-\frac{1}{2}ca+3c_{0}(a-1)
+\frac{3}{2}(a^{2}-1)(r-1)}
((c/2)(2\tau)+d)^{4}\eta(4\tau)^{8},
\end{align*}
and thus,
\begin{align*}
    &\eta\Big(\frac{a\tau+b}{c\tau+d}\Big)^8
    \eta\Big(\frac{a(2\tau)+2b}{(c/2)(2\tau)+d}\Big)^8\\
    &=\omega^{3ab+\frac{3}{2}cd(1-a^{2})+6c_{0}(a-1)+\frac{3}{2}(2r-1)(a^{2}-1)}(c\tau+d)^{8}\eta(\tau)^{8}\eta(2\tau)^{8}=(c\tau+d)^{8}\eta(\tau)^{8}\eta(2\tau)^{8}.
\end{align*}
Combining these implies that
$$
  u\Big(\frac{a\tau+b}{c\tau+d}\Big)
  =\Big(\frac{\eta(\frac{a\tau+b}{c\tau+d})^{2}
  \eta(\frac{a(2\tau)+2b}{(c/2)(2\tau)+d})^{2}}{8E_{2}
  (\frac{a(2\tau)+b}{(c/2)(2\tau)+d})-4E_{2}(\frac{a\tau+b}{c\tau+d})}
  \Big)^{4}=u(\tau)
\q\t{for $\begin{pmatrix}
    a&b\\c&d
\end{pmatrix}\in\Gamma_{0}(2)$.}$$
 By (4.3),
$u(\tau)=u(-\f 1{2\tau})$.
Therefore, $u(\tau)$ is invariant under $\Gamma_{0}(2)^{+}$. Clearly, $u(\tau)=q+O(q^{2})$ and thus has a unique zero at the cusp $[i\infty]$ of order~1, where $q=e^{2\pi i\tau}$. So by the Residue Theorem for compact Riemann surfaces, $u(\tau)$ must be a Hauptmodul for $\Gamma_{0}(2)^{+}$. In addition, by the Riemann--Roch theorem, one can show that $2E_{2}(\tau)-E_{2}(\tau)$ has a zero of order~$\frac{1}{2}$ at the elliptic point $\frac{1+i}{2}$ of $\Gamma_{0}(2)$, which is of period~2, and so, with respect to $\Gamma_{0}(2)^{+}$, it has a zero of order~$\frac{1}{4}$ at $\frac{1+i}{2}$, which is an elliptic point of period~4 of $\Gamma_{0}(2)^{+}$. The proof is now complete.

\v2
\par{\bf Lemma 4.2} {\it Let $u(\tau)$ be given in $(4.1)$. Then
\begin{align*}&u\Ls{\sqrt{-2}}2=\f 1{256},\q u\Ls{1+\sqrt{-3}}2=-\f 1{144},\q u\Ls{1+i}4=\f 1{648},
\\& u\Ls{1+\sqrt{-7}}2=-\f 1{3969},\q u\Ls{1+\sqrt{-7}}4=\f 1{81}.\end{align*}}
\par{\it Proof.} Since $f_2(-\f 1{\tau})=f_1(\tau)$ and $f_1(2\tau)f_2(\tau)=\sqrt 2$, we see that
$f_2\sls{\sqrt{-2}}2=f_1(\sqrt {-2})$
and so $f_2\sls{\sqrt{-2}}2^2=f_1(\sqrt {-2})f_2\sls{\sqrt{-2}}2=\sqrt 2$. Thus, $f_2\sls{\sqrt{-2}}2^{24}=(\sqrt 2)^{12}=64$ and so
$$u\Ls{\sqrt{-2}}2=\f{f_2\sls{\sqrt{-2}}2^{24}}{(f_2\sls{\sqrt{-2}}2^{24}+64)^2}
=\f{64}{(64+64)^2}=\f 1{256}.$$
\par By [7, p.261], $\gamma_2(\f{1+\sqrt{-3}}2)=0$. Thus
$f_2(\f{1+\sqrt{-3}}2)^{24}+16=\gamma_2(\f{1+\sqrt{-3}}2)
f_2(\f{1+\sqrt{-3}}2)^8=0$.
Hence $f_2(\f{1+\sqrt{-3}}2)^{24}=-16$ and so
$u\ls{1+\sqrt{-3}}2=\f{-16}{(-16+64)^2}=-\f 1{144}$ by (4.1).
\par From [7, p.261], $\gamma_2(i)=12$. Thus,
$$\gamma_2\Ls{-1+i}2=\gamma_2\Big(-\f 1{\f{-1+i}2}\Big)=\gamma_2(1+i)=\omega^2\gamma_2(i)=12\omega^2.$$
By Lemma 2.2, $-f\ls{-1+i}2^8,f_1\ls{-1+i}2^8$ and $f_2\ls{-1+i}2^8$ are the three roots of $x^3-12\omega^2x+16=0$. Notice that
$$x^3-12\omega^2x+16=\Ls x{\omega}^3-12\f x{\omega}+16=
\Big(\f x{\omega}-2\Big)^2\Big(\f x{\omega}+4\Big).$$
 We have
$$-f\Ls{-1+i}2^8,\; f_1\Ls{-1+i}2^8,\;f_2\Ls{-1+i}2^8\in\{2\omega,-4\omega\}.$$
Since $f_2(-\f 1{\tau})=f_1(\tau)$, we see that
$f_2(\f{-1+i}2)=f_1(1+i)$ and so
$f_2(\f{-1+i}2)f_2(\f{1+i}2)=f_1(1+i)f_2(\f{1+i}2)=\sqrt 2$.
This together with Lemma 2.1 yields $|f_2\sls{-1+i}2|^2=f_2\sls{-1+i}2\overline{f_2\sls{-1+i}2}$ $=f_2\sls{-1+i}2
f_2\sls{1+i}2=\sqrt 2$. Therefore, $|f_2\sls{-1+i}2^8|=(\sqrt 2)^4=4$ and so $f_2\sls{-1+i}2^8=-4\omega$. Since $-4\omega$ is a single root of $x^3-12\omega^2x+16=0$, we must have
$f_1\sls{-1+i}2^8=2\omega$. Applying Lemma 2.1 gives
$f_1\ls{1+i}2^8=\overline{f_1\ls{-1+i}2}^8=\overline{f_1\ls{-1+i}2^8}=
\overline{2\omega}=2\omega^2$
and so $f_1\ls{1+i}2^{24}=8$. Hence
$u\ls{1+i}4=\f{8}{(8+64)^2}=\f 1{648}$ by (4.2).

\par By [7, (12.20)], $\gamma_2\big(\f{3+\sqrt{-7}}2\big)=-15$. Thus,
$\gamma_2(\f{1+\sqrt{-7}}2)=\omega\gamma_2\big(1+\f{1+\sqrt{-7}}2\big)
=-15\omega$. In view of Lemma 2.2, $-f\ls{1+\sqrt{-7}}2^8,\ f_1\ls{1+\sqrt{-7}}2^8$ and
$f_2\ls{1+\sqrt{-7}}2^8$ are the distinct roots of $x^3+15\omega x+16=
0$. Observe that
$$x^3+15\omega x+16=\Ls x{\omega^2}^3+15\f x{\omega^2}+16
=\Big(\f x{\omega^2}+1\Big)\Big(\big(\f x{\omega^2}-\f 12\big)^2+\f{63}4\Big).$$
 We then obtain
$$-f\ls{1+\sqrt{-7}}2^8, f_1\ls{1+\sqrt{-7}}2^8, f_2\big(\f{1+\sqrt{-7}}2\big)^8
\in\Big\{-\omega^2,\f{1+3\sqrt{-7}}2\omega^2,
\f{1-3\sqrt{-7}}2\omega^2\Big\}$$
and so
$$ f_1\ls{1+\sqrt{-7}}2^{24},\ f_2\big(\f{1+\sqrt{-7}}2\big)^{24}
\in\Big\{-1,-\f{47+45\sqrt{-7}}2,
-\f{47-45\sqrt{-7}}2\Big\}.$$
This implies that $f_2\sls{1+\sqrt{-7}}2^{24}$ is real if and  only if $f_2\sls{1+\sqrt{-7}}2^{24}=-1$.
Since $e^{2\pi i\f{1+\sqrt{-7}}2}=e^{\pi i-\sqrt 7\pi}=-e^{-\sqrt 7\pi}$ we see that
$$f_2\Ls{1+\sqrt{-7}}2^{24}=2^{12}\big(-e^{-\sqrt 7\pi}\big)
\prod_{n=1}^{\infty}\Big(1+\big(-e^{-\sqrt 7\pi}\big)^n\Big)^{24}$$
is real. Therefore,
$$f_2\Ls{1+\sqrt{-7}}2^{24}=-1\q \t{and so}\q
f_2\Ls{1+\sqrt{-7}}2^8=-\omega^2.$$
It then follows that
$$u\Ls{1+\sqrt{-7}}2=\f{f_2\ls{1+\sqrt{-7}}2^{24}}{\big(f_2\sls{1+\sqrt{-7}}2^{24}
+64\big)^2}=\f{-1}{(64-1)^2}=-\f 1{3969}.$$
\par Since the roots of $x^3+15\omega x+16=
0$ are distinct, from the above we deduce that
$$f_1\Ls{1+\sqrt{-7}}2^8=\f{1\pm 3\sqrt{-7}}2\omega^2
\q\t{and so}\q
f_1\Ls{1+\sqrt{-7}}2^{24}=-\f{47\pm 45\sqrt{-7}}2.$$
Hence
$$f_2\Ls{1+\sqrt{-7}}4^{24}=\f{2^{12}}{f_1\ls{1+\sqrt{-7}}2^{24}}
=\f{2^{12}}{-\f{47\pm 45\sqrt{-7}}2}=-\f{47\mp 45\sqrt {-7}}2$$
and so
\begin{align*}u\Ls{1+\sqrt{-7}}4&=\f{f_2\sls{1+\sqrt{-7}}4^{24}}
{\big(f_2\sls{1+\sqrt{-7}}4^{24}
+64\big)^2}=\f 1{f_2\sls{1+\sqrt{-7}}4^{24}+128+\f{2^{12}}{f_2\sls{1+\sqrt{-7}}4^{24}}}
\\&=\f 1{-\f{47\mp 45\sqrt{-7}}2+128-\f{47\pm 45\sqrt{-7}}2}=\f 1{81}.
\end{align*}
The proof is now complete.
\v2

\par{\bf Remark 4.1} Lemma 4.2 was stated by Beukers in [2, p.30] without proof.
\v2
\par Let $u(\tau)$ be given in $(4.1)$.
By Lemma 4.1, both
$u(\tau)$ and $\sum_{n=0}^{\infty}\binom{2n}{n}^2\b{4n}{2n}u(\tau)^{n}
$ satisfy the assumptions in Beukers' theorem.
\v2
\par{\bf Proof of Theorem 1.5.} Since $p=x^2+2y^2=(x+y\sqrt{-2})(x-y\sqrt{-2})$ for $x,y\in\Bbb Z$, one may choose the sign of $y$ so that $x+y\sqrt{-2}\notin p\Bbb Z_p$.
 Set
 $$a_n=\b{2n}n^2\b{4n}{2n},\q c=2y,\q d=x,\q \alpha=\f{\sqrt{-2}}2\q\t{and}\q N=2.$$
 Then clearly
 $(c,d)=1,\ N\mid c,$ $c\alpha\bar{\alpha},\ c(\alpha+\bar{\alpha})\in\Bbb Z$,
$p=(c\alpha+d)(c\bar{\alpha}+d)$ and $c\alpha+d\notin p\Bbb Z_p$.
 By Lemma 4.2, $u\sls{\sqrt{-2}}2=\f 1{256}$.
 Hence, applying Beukers' theorem(i) and Lemma 2.3 we obtain
 \begin{align*}\sum_{n=0}^{p-1}\f{\b{2n}n^2\b{4n}{2n}}{256^n}
 &=\sum_{n=0}^{p-1}\b{2n}n^2\b{4n}{2n}u\Ls{\sqrt{-2}}2^n\e \Big(2y\cdot \f{\sqrt{-2}}2+x\Big)^2
 \\&\e 4x^2-2p-\f{p^2}{4x^2}\mod {p^3}.\end{align*}

\par{\bf Proof of Theorem 1.6.} Since $p=x^2+3y^2=(x+y\sqrt{-3})(x-y\sqrt{-3})$ for $x,y\in\Bbb Z$, one may choose the sign of $y$ so that $x+y\sqrt{-3}\notin p\Bbb Z_p$.
 Set
 $$a_n=\b{2n}n^2\b{4n}{2n},\q c=2y,\q d=x-y,\q \alpha=\f{1+\sqrt{-3}}2\q\t{and}\q N=2.$$
 Then clearly
 $(c,d)=1,\ N\mid c,$ $c\alpha\bar{\alpha},\ c(\alpha+\bar{\alpha})\in\Bbb Z$,
$p=(c\alpha+d)(c\bar{\alpha}+d)$ and $c\alpha+d\notin p\Bbb Z_p$.
 By Lemma 4.2, $u\sls{1+\sqrt{-3}}2=-\f 1{144}$.
 Hence, applying Beukers' theorem(i) and Lemma 2.3 we obtain
 \begin{align*}\sum_{n=0}^{p-1}\f{\b{2n}n^2\b{4n}{2n}}{(-144)^n}
 &=\sum_{n=0}^{p-1}\b{2n}n^2\b{4n}{2n}u\Ls{1+\sqrt{-3}}2^n\e \Big(2y\cdot \f{1+\sqrt{-3}}2+x-y\Big)^2
 \\&=(x+y\sqrt {-3})^2\e 4x^2-2p-\f{p^2}{4x^2}\mod {p^3}.\end{align*}

 \par{\bf Proof of Theorem 1.7.} Since $p=x^2+4y^2=(x+2yi)(x-2yi)$ for $x,y\in\Bbb Z$, one may choose the sign of $y$ so that $x+2yi\notin p\Bbb Z_p$.  Set
 $$a_n=\b{2n}n^2\b{4n}{2n},\q c=8y,\q d=x-2y,\q \alpha=\f{1+i}4\q\t{and}\q N=2.$$
 Then clearly
 $(c,d)=1,\ N\mid c,$ $c\alpha\bar{\alpha},\ c(\alpha+\bar{\alpha})\in\Bbb Z$,
$p=(c\alpha+d)(c\bar{\alpha}+d)$ and $c\alpha+d\notin p\Bbb Z_p$.
 By Lemma 4.2, $u\sls{1+i}4=\f 1{648}$.
 Hence, applying Beukers' theorem(i) and Lemma 2.3 we obtain
 \begin{align*}\sum_{n=0}^{p-1}\f{\b{2n}n^2\b{4n}{2n}}{648^n}
 &=\sum_{n=0}^{p-1}\b{2n}n^2\b{4n}{2n}u\Ls{1+i}4^n\e \Big(8y\cdot \f{1+i}4+x-2y\Big)^2
 \\&=(x+y\sqrt{-4})^2\e 4x^2-2p-\f{p^2}{4x^2}\mod {p^3}.\end{align*}

 \par{\bf Proof of Theorem 1.8.} Since $p=x^2+7y^2=(x+y\sqrt{-7})(x-y\sqrt{-7})$ for $x,y\in\Bbb Z$, one may choose the sign of $y$ so that $x+y\sqrt{-7}\notin p\Bbb Z_p$.
  Set
 $$a_n=\b{2n}n^2\b{4n}{2n},\q c=2y,\q d=x-y,\q \alpha=\f{1+\sqrt{-7}}2\q\t{and}\q N=2.$$
 Then clearly
 $(c,d)=1,\ N\mid c,$ $c\alpha\bar{\alpha},\ c(\alpha+\bar{\alpha})\in\Bbb Z$,
$p=(c\alpha+d)(c\bar{\alpha}+d)$ and $c\alpha+d\notin p\Bbb Z_p$.
 By Lemma 4.2, $u\sls{1+\sqrt{-7}}2=-\f 1{3969}$.
 Hence, applying Beukers' theorem(i) and Lemma 2.3 we obtain
 \begin{align*}\sum_{n=0}^{p-1}\f{\b{2n}n^2\b{4n}{2n}}{(-3969)^n}
 &=\sum_{n=0}^{p-1}\b{2n}n^2\b{4n}{2n}u\Ls{1+\sqrt{-7}}2^n\e \Big(2y\cdot \f{1+\sqrt{-7}}2+x-y\Big)^2
 \\&=(x+y\sqrt {-7})^2\e 4x^2-2p-\f{p^2}{4x^2}\mod {p^3}.\end{align*}
 On the other hand, setting
 $$a_n=\b{2n}n^2\b{4n}{2n},\q c=4y,\q d=x-y,\q \alpha=\f{1+\sqrt{-7}}4\q\t{and}\q N=2$$
 one finds that
 $(c,d)=1,\ N\mid c,$ $c\alpha\bar{\alpha},\ c(\alpha+\bar{\alpha})\in\Bbb Z$,
$p=(c\alpha+d)(c\bar{\alpha}+d)$ and $c\alpha+d\notin p\Bbb Z_p$.
 By Lemma 4.2, $u\sls{1+\sqrt{-7}}4=\f 1{81}$.
 Hence, applying Beukers' theorem(i) and Lemma 2.3 we obtain
 \begin{align*}\sum_{n=0}^{p-1}\f{\b{2n}n^2\b{4n}{2n}}{81^n}
 &=\sum_{n=0}^{p-1}\b{2n}n^2\b{4n}{2n}u\Ls{1+\sqrt{-7}}4^n\e \Big(4y\cdot \f{1+\sqrt{-7}}4+x-y\Big)^2
 \\&=(x+y\sqrt {-7})^2\e 4x^2-2p-\f{p^2}{4x^2}\mod {p^3}.\end{align*}

\par{\bf Lemma 4.3} {\it We have
\begin{align*}&u\Ls{1+3i}2=-\f 1{12288},\q u\Ls{1+5i}2=-\f 1{6635520},
\q u\Ls{1+\sqrt{-5}}2=-\f 1{1024},\\& u\Ls{1+\sqrt{-13}}2=-\f 1{82944} ,\q u\Ls{1+\sqrt{-37}}2=-\f 1{14112^2},
\\&u\Ls{\sqrt{-6}}2=\f 1{48^2},\q u\Ls{\sqrt{-10}}2=\f 1{12^4},\q u\Ls{3\sqrt{-2}}2=\f 1{28^4},\\& u\Ls{\sqrt{-22}}2=\f 1{1584^2},
\q u\Ls{\sqrt{-58}}2=\f 1{396^4}.\end{align*}}
\par{\it Proof.} For positive integers $n$ let
$G_n=\f 1{2^{\f 14}}f(\sqrt{-n})$. Since $f(\tau+1)=\zeta_{48}^{-1}f(\tau)$, we have
$$f_1(1+\sqrt{-n})^{24}=-f(\sqrt{-n})^{24}=-(2^{\f 14}G_n)^{24}=-64G_n^{24}$$
and so
\begin{align*}u\Ls{1+\sqrt{-n}}2
&=\f{f_1(1+\sqrt{-n})^{24}}{(f_1(1+\sqrt{-n})^{24}+64)^2}
=\f 1{f_1(1+\sqrt{-n})^{24}+128+\f{4096}{f_1(1+\sqrt{-n})^{24}}}
\\&=\f 1{-64G_n^{24}+128+\f{4096}{-64G_n^{24}}}=-\f 1{64(G_n^{12}-G_n^{-12})^2}.
\end{align*}
From [1, pp.189-191] we know that
\begin{align*}&G_5=\Ls{1+\sqrt 5}2^{\f 14},\q G_9=\Ls{1+\sqrt 3}{\sqrt 2}^{\f 13},\q G_{13}=\Ls{3+\sqrt{13}}2^{\f 14},
\\& G_{25}=\f{1+\sqrt 5}2,\q
G_{37}=(6+\sqrt{37})^{\f 14}.\end{align*}
Thus,
\begin{align*}
\\ u\Ls{1+\sqrt{-5}}2&=-\f 1{64(G_{5}^{12}-G_{5}^{-12})^2}=-\f 1{64(\sls{1+\sqrt 5}2^3-\sls{1+\sqrt 5}2^{-3})^2}
\\&=-\f 1{64(2+\sqrt 5+(2-\sqrt 5))^2}=-\f 1{1024},\end{align*}
\begin{align*}
\\u\Ls{1+3i}2&=-\f 1{64(G_9^{12}-G_9^{-12})^2}=-\f 1{64(\f{(1+\sqrt 3)^4}4-\f 4{(1+\sqrt 3)^4})^2}
\\&=-\f 1{64(7+4\sqrt 3-(7-4\sqrt 3))^2}=-\f 1{12288},
 \\ u\Ls{1+5i}2&=-\f 1{64(G_{25}^{12}-G_{25}^{-12})^2}=-\f 1{64(\sls{1+\sqrt 5}2^{12}-\sls{1+\sqrt 5}2^{-12})^2}
\\&=-\f 1{64(\sls{1+\sqrt 5}2^{12}-\sls{1-\sqrt 5}2^{12})^2}
=-\f 1{64\cdot 5\cdot 144^2}=-\f 1{6635520},
\\u\Ls{1+\sqrt{-13}}2&=-\f 1{64(G_{13}^{12}-G_{13}^{-12})^2}=-\f 1{64(\f{(3+\sqrt {13})^3}8-\f 8{(3+\sqrt {13})^3})^2}
\\&=-\f 1{64(18+5\sqrt{13}+18-5\sqrt{13})^2}=-\f 1{82944},
\\u\Ls{1+\sqrt{-37}}2&=-\f 1{64(G_{37}^{12}-G_{37}^{-12})^2}=-\f 1{64((6+\sqrt {37})^3-(6+\sqrt{37})^{-3})^2}
\\&=-\f 1{64(882+145\sqrt{37}+(882-145\sqrt{37}))^2}=-\f 1{14112^2}.
\end{align*}
\par For positive integers $n$ let
$g_n=\f 1{2^{\f 14}}f_1(\sqrt {-n}).$ Then $f_1(\sqrt{-n})^{24}=64g_n^{24}$. Hence
\begin{align*}u\Ls{\sqrt{-n}}2=\f{f_1(\sqrt{-n})^{24}}{(f_1(\sqrt{-n})^{24}+64)^2}
=\f{64g_n^{24}}{(64g_n^{24}+64)^2}
=\f 1{64 (g_n^{12}+g_n^{-12})^2}.\end{align*}
From [1, pp.200-201] we know that
\begin{align*}&g_6=(1+\sqrt 2)^{\f 16},\q g_{10}=\Ls{1+\sqrt 5}2^{\f 12},
\q g_{18}=(\sqrt 2+\sqrt 3)^{\f 13},
\\& g_{22}=(1+\sqrt 2)^{\f 12},
\q g_{58}=\Ls{5+\sqrt{29}}2^{\f 12}.\end{align*}
Thus,
\begin{align*}&g_6^{12}=3+2\sqrt 2,\q  g_{10}^{12}=\Ls{1+\sqrt 5}2^6=
9+4\sqrt 5,\q
g_{18}^{12}=(\sqrt 2+\sqrt 3)^4
=49+20\sqrt 6,\\& g_{22}^{12}=(1+\sqrt 2)^6=99+70\sqrt 2,\q
g_{58}^{12}=\Ls{5+\sqrt{29}}2^6=9801+1820\sqrt{29}.\end{align*}
Hence,
\begin{align*}
&u\Ls{\sqrt{-6}}2=\f 1{64(g_6^{12}+g_6^{-12})^2}=\f 1{64(3+2\sqrt 2+3-2\sqrt 2)^2}=\f 1{48^2},
\\&u\Ls{\sqrt{-10}}2=\f 1{64(g_{10}^{12}+g_{10}^{-12})^2}=\f 1{64(9+4\sqrt 5+9-4\sqrt 5)^2}=\f 1{12^4},
\\&u\Ls{3\sqrt{-2}}2=\f 1{64(g_{18}^{12}+g_{18}^{-12})^2}=\f 1{64(49+20\sqrt {6}+49-20\sqrt {6})^2}=\f 1{28^4},
\\&u\Ls{\sqrt{-22}}2=\f 1{64(g_{22}^{12}+g_{22}^{-12})^2}=\f 1{64(99+70\sqrt 2+99-70\sqrt 2)^2}=\f 1{1584^2},
\\&u\Ls{\sqrt{-58}}2=\f 1{64(g_{58}^{12}+g_{58}^{-12})^2}
=\f 1{64(9801+1820\sqrt{29}+9801-1820\sqrt{29})^2}
 =\f 1{396^4}.
\end{align*}
\par{\bf Remark 4.2} Lemma 4.3 was stated by Beukers in [2, p.30] without proof.
\v2
\par Let $u(\tau)$ be given in $(4.1)$.
By Lemma 4.1, both
$u(\tau)$ and $\sum_{n=0}^{\infty}\binom{2n}{n}^2\b{4n}{2n}u(\tau)^{n}
$ satisfy the assumptions in Beukers' theorem.
\v2
\par{\bf Proof of Theorem 1.9.} Since $p=x^2+2y^2=(x+y\sqrt{-2})(x-y\sqrt{-2})$ for $x,y\in\Bbb Z$, one may choose the sign of $y$ so that $x+y\sqrt{-2}\notin p\Bbb Z_p$.
For $p\e 1\mod 3$ we have $3\mid y$.  Set
 $$a_n=\b{2n}n^2\b{4n}{2n},\q c=\f{2y}3,\q d=x,\q \alpha=\f{3\sqrt{-2}}2\q\t{and}\q N=2.$$
 Then clearly
 $(c,d)=1,\ N\mid c,$ $c\alpha\bar{\alpha},\ c(\alpha+\bar{\alpha})\in\Bbb Z$,
$p=(c\alpha+d)(c\bar{\alpha}+d)$ and $c\alpha+d\notin p\Bbb Z_p$.
 By Lemma 4.3, $u\sls{3\sqrt{-2}}2=\f 1{28^4}$.
 Hence, applying Beukers' theorem(i) and Lemma 2.3 we obtain
 \begin{align*}\sum_{n=0}^{p-1}\f{\b{2n}n^2\b{4n}{2n}}{28^{4n}}
 &=\sum_{n=0}^{p-1}\b{2n}n^2\b{4n}{2n}u\Ls{3\sqrt{-2}}2^n\e \Big(\f{2y}3\cdot \f{3\sqrt{-2}}2+x\Big)^2
 \\&\e 4x^2-2p-\f{p^2}{4x^2}\mod {p^3}.\end{align*}
\par Now assume $p\e 2\mod 3$. Then $3\mid x$. Setting
 $$a_n=\b{2n}n^2\b{4n}{2n},\q c=\f{x}3,\q d=-y,\q \alpha=\f{3\sqrt{-2}}2\q\t{and}\q N=2.$$
we find $(c,dN)=1,\ Nc\alpha\bar{\alpha}\in\Bbb Z,$  $c(\alpha+\bar{\alpha})\in\Bbb Z$,
$p=N(c\alpha+d)(c\bar{\alpha}+d)$ and $c\alpha+d\notin p\Bbb Z_p$.
Since $u\ls{3\sqrt{-2}}2=\f 1{28^4}$ by Lemma 4.3, using Beukers' theorem(ii) and Lemma 2.3 we deduce that
\begin{align*}\sum_{n=0}^{p-1}\f{\b{2n}n^2\b{4n}{2n}}{28^{4n}}
 &=\sum_{n=0}^{p-1}\b{2n}n^2\b{4n}{2n}u\Ls{3\sqrt{-2}}2^n\e -2\Big(\f{x}3\cdot \f{3\sqrt{-2}}2-y\Big)^2
 \\&\e (x+\sqrt{-2}y)^2\e 4x^2-2p-\f{p^2}{4x^2}\mod {p^3}.\end{align*}

\par{\bf Proofs of Theorems 1.10-1.11.} Let $m\in\{5,9,13,25,37\}$ and
$D(m)=-1024,-12288,-82944,$ $-6635520,$ $-14112^2$ according as $m=5,9,13,25,37$. For $p=x^2+my^2$ we have $p=(x+y\sqrt{-m})(x-y\sqrt{-m})$ for $x,y\in\Bbb Z$, one may choose the sign of $y$ so that $x+y\sqrt{-m}\notin p\Bbb Z_p$.
Set
 $$a_n=\b{2n}n^2\b{4n}{2n},\q c=2y,\q d=x-y,\q \alpha=\f{1+\sqrt{-m}}2\q\t{and}\q N=2.$$
 Then clearly
 $(c,d)=1,\ N\mid c,$ $c\alpha\bar{\alpha},\ c(\alpha+\bar{\alpha})\in\Bbb Z$,
$p=(c\alpha+d)(c\bar{\alpha}+d)$ and $c\alpha+d\notin p\Bbb Z_p$.
 By Lemma 4.3, $u\sls{1+\sqrt{-m}}2=\f 1{D(m)}$.
 Hence, applying Beukers' theorem(i) and Lemma 2.3 we obtain
 \begin{align*}\sum_{n=0}^{p-1}\f{\b{2n}n^2\b{4n}{2n}}{D(m)^n} &=\sum_{n=0}^{p-1}\b{2n}n^2\b{4n}{2n}u\Ls{1+\sqrt{-m}}2^n\e \Big(2y\cdot \f{1+\sqrt{-m}}2+x-y\Big)^2
 \\&=(x+\sqrt{-m}y)^2\e 4x^2-2p-\f{p^2}{4x^2}\mod {p^3}.\end{align*}
\par For $2p=x^2+my^2$ we have $2\nmid xy$. One may choose the sign of $y$ so that $x+y\sqrt{-m}\notin p\Bbb Z_p$.
Setting
 $$a_n=\b{2n}n^2\b{4n}{2n},\q c=y,\q d=\f{x-y}2,\q \alpha=\f{1+\sqrt{-m}}2\q\t{and}\q N=2.$$
we find $(c,dN)=1,\ Nc\alpha\bar{\alpha}\in\Bbb Z,$  $c(\alpha+\bar{\alpha})\in\Bbb Z$,
$p=N(c\alpha+d)(c\bar{\alpha}+d)$ and $c\alpha+d\notin p\Bbb Z_p$.
Since $u\ls{1+\sqrt{-m}}2=\f 1{D(m)}$ by Lemma 4.3, using Beukers' theorem(ii) and Lemma 2.3 we deduce that
\begin{align*}\sum_{n=0}^{p-1}\f{\b{2n}n^2\b{4n}{2n}}{D(m)^n} &=\sum_{n=0}^{p-1}\b{2n}n^2\b{4n}{2n}u\Ls{1+\sqrt{-m}}2^n\e -2\Big(y\cdot \f{1+\sqrt{-m}}2+\f{x-y}2\Big)^2
 \\&=-\f 12(x+\sqrt{-m}y)^2\e -\f 12\Big(4x^2-4p-\f{p^2}{x^2}\Big)\\&=-2x^2+2p+\f{p^2}{2x^2}\mod {p^3}.\end{align*}

\par{\bf Proof of Theorem 1.12.}  For $p=x^2+2my^2$ we have $p=(x+y\sqrt{-2m})(x-y\sqrt{-2m})$ for $x,y\in\Bbb Z$, one may choose the sign of $y$ so that $x+y\sqrt{-2m}\notin p\Bbb Z_p$.
Set
 $$a_n=\b{2n}n^2\b{4n}{2n},\q c=2y,\q d=x,\q \alpha=\f{\sqrt{-2m}}2\q\t{and}\q N=2.$$
 Then clearly
 $(c,d)=1,\ N\mid c,$ $c\alpha\bar{\alpha},\ c(\alpha+\bar{\alpha})\in\Bbb Z$,
$p=(c\alpha+d)(c\bar{\alpha}+d)$ and $c\alpha+d\notin p\Bbb Z_p$.
 By Lemma 4.3, $u\sls{\sqrt{-2m}}2=\f 1{F(m)}$.
 Hence, applying Beukers' theorem(i) and Lemma 2.3 we obtain
 \begin{align*}\sum_{n=0}^{p-1}\f{\b{2n}n^2\b{4n}{2n}}{F(m)^n} &=\sum_{n=0}^{p-1}\b{2n}n^2\b{4n}{2n}u\Ls{\sqrt{-2m}}2^n\e \Big(2y\cdot \f{\sqrt{-2m}}2+x\Big)^2
 \\&=(x+\sqrt{-2m}y)^2\e 4x^2-2p-\f{p^2}{4x^2}\mod {p^3}.\end{align*}
\par For $p=2x^2+my^2$ we have $2\nmid y$ and $2p=(2x)^2+2my^2=(2x+y\sqrt{-2m})(2x-y\sqrt{-2m})$. One may choose the sign of $y$ so that $2x+y\sqrt{-2m}\notin p\Bbb Z_p$.
Setting
 $$a_n=\b{2n}n^2\b{4n}{2n},\q c=y,\q d=x,\q \alpha=\f{\sqrt{-2m}}2\q\t{and}\q N=2.$$
we find $(c,dN)=1,\ Nc\alpha\bar{\alpha}\in\Bbb Z,$  $c(\alpha+\bar{\alpha})\in\Bbb Z$,
$p=N(c\alpha+d)(c\bar{\alpha}+d)$ and $c\alpha+d\notin p\Bbb Z_p$.
Since $u\ls{\sqrt{-2m}}2=\f 1{F(m)}$ by Lemma 4.3, using Beukers' theorem(ii) and Lemma 2.3 we deduce that
\begin{align*}\sum_{n=0}^{p-1}\f{\b{2n}n^2\b{4n}{2n}}{F(m)^n} &=\sum_{n=0}^{p-1}\b{2n}n^2\b{4n}{2n}u\Ls{\sqrt{-2m}}2^n\e -2\Big(y\cdot \f{\sqrt{-2m}}2+x\Big)^2
 \\&=-\f 12(2x+\sqrt{-2m}y)^2\e -\f 12\Big(4(2x)^2-4p-\f{4p^2}{4(2x)^2}\Big)\\&=-8x^2+2p+\f{p^2}{8x^2}\mod {p^3}.\end{align*}

\section*{5. Proofs of congruences involving $\b{2k}k\b{3k}k\b{6k}{3k}$}

\par\q For given positive integer $k$ let $(a)_k=a(a+1)\cdots(a+k-1)$. Then $(a)_k=(-1)^k\b{-a}kk!$. From [19] we know that
$$\f{\ls 12_k\ls 16_k \ls 56_k}{(1)_k^3}
=(-1)^k\b {-\f 12}k\b {-\f 16}k\b{-\f 56}k=\f{\b{2k}k}{4^k}
\cdot\f{\b{3k}k\b{6k}{3k}}{432^k}=\f{\b{2k}k\b{3k}k\b{6k}{3k}}{12^{3k}}.$$

\par{\bf Proof of Theorem 1.13.} Suppose that $p=4n+1=x^2+4y^2$ so that $x+y\sqrt{-4}\not\in p\Bbb Z_p$. Taking $D=4$ in [2, Theorem 1.27] and $d_K=-4$ in [7, (12.20)] gives
$$\sum_{k=0}^{p-1}\f{\b{2k}k\b{3k}k\b{6k}{3k}}{12^{3k}}\e \pm (x+y\sqrt{-4})^2\mod {p^3}.$$
Applying Lemma 2.3 we get
$$\sum_{k=0}^{p-1}\f{\b{2k}k\b{3k}k\b{6k}{3k}}{12^{3k}}\e
\pm \Big(4x^2-2p-\f{p^2}{4x^2}\Big)\mod{p^3}.$$
Since Mortenson[13] proved that $\sum_{k=0}^{p-1}\f{\b{2k}k\b{3k}k\b{6k}{3k}}{12^{3k}}\e \ls {-3}p4x^2\mod {p^2}$, we must have
$$\sum_{k=0}^{p-1}\f{\b{2k}k\b{3k}k\b{6k}{3k}}{12^{3k}}\e
\Ls {-3}p \Big(4x^2-2p-\f{p^2}{4x^2}\Big)\mod{p^3}.$$
\par On the other hand, taking $D=16$ in [2, Theorem 1.27] and $d_K=-16$ in [7, (12.20)] gives
$$\sum_{k=0}^{p-1}\f{\b{2k}k\b{3k}k\b{6k}{3k}}{66^{3k}}\e \pm (x+y\sqrt{-4})^2\mod {p^3}.$$
 Applying Lemma 2.3 we get
$$\sum_{k=0}^{p-1}\f{\b{2k}k\b{3k}k\b{6k}{3k}}{66^{3k}}\e
\pm \Big(4x^2-2p-\f{p^2}{4x^2}\Big)\mod{p^3}.$$
By [19, Theorem 4.3],
$\sum_{k=0}^{p-1}\f{\b{2k}k\b{3k}k\b{6k}{3k}}{66^{3k}}
\e\ls {33}p4x^2\mod p.$
Hence,
$$\sum_{k=0}^{p-1}\f{\b{2k}k\b{3k}k\b{6k}{3k}}{66^{3k}}
\e\Ls {33}p\Big(4x^2-2p-\f{p^2}{4x^2}\Big)\mod{p^3}.$$

\par{\bf Proof of Theorem 1.14.} Suppose that $p=3n+1=x^2+3y^2(x,y\in\Bbb Z)$ and we choose the sign of $y$ so that $x+y\sqrt{-3}\not\in p\Bbb Z_p$. Taking $D=12$ in [2, Theorem 1.27] and $d_K=-12$ in [7, (12.20)] gives
$$\sum_{k=0}^{p-1}\f{\b{2k}k\b{3k}k\b{6k}{3k}}{54000^k}\e \pm (x+y\sqrt{-3})^2\mod {p^3}.$$
 Applying Lemma 2.3 we get
$$\sum_{k=0}^{p-1}\f{\b{2k}k\b{3k}k\b{6k}{3k}}{54000^k}\e
\pm \Big(4x^2-2p-\f{p^2}{4x^2}\Big)\mod{p^3}.$$
By [19, Theorem 4.5],
$\sum_{k=0}^{p-1}\f{\b{2k}k\b{3k}k\b{6k}{3k}}{54000^k}
\e\ls {5}p4x^2\mod p.$
Hence,
$$\sum_{k=0}^{p-1}\f{\b{2k}k\b{3k}k\b{6k}{3k}}{54000^k}
\e\Ls {5}p\Big(4x^2-2p-\f{p^2}{4x^2}\Big)\mod{p^3}.$$

\par{\bf Proof of Theorem 1.15.} Suppose that $p\e 1,3\mod 8$ and so $p=x^2+2y^2(x,y\in\Bbb Z)$. We choose the sign of $y$ so that $x+y\sqrt{-2}\not\in p\Bbb Z_p$. Taking $D=8$ in [2, Theorem 1.27] and $d_K=-8$ in [7, (12.20)] gives
$$\sum_{k=0}^{p-1}\f{\b{2k}k\b{3k}k\b{6k}{3k}}{20^{3k}}\e \pm (x+y\sqrt{-2})^2\mod {p^3}.$$
 Applying Lemma 2.3 we get
$$\sum_{k=0}^{p-1}\f{\b{2k}k\b{3k}k\b{6k}{3k}}{20^{3k}}\e
\pm \Big(4x^2-2p-\f{p^2}{4x^2}\Big)\mod{p^3}.$$
By [19, Theorem 4.4],
$\sum_{k=0}^{p-1}\f{\b{2k}k\b{3k}k\b{6k}{3k}}{20^{3k}}
\e\ls {-5}p4x^2\mod p.$
Hence,
$$\sum_{k=0}^{p-1}\f{\b{2k}k\b{3k}k\b{6k}{3k}}{20^{3k}}
\e\Ls {-5}p\Big(4x^2-2p-\f{p^2}{4x^2}\Big)\mod{p^3}.$$
\v2
\par{\bf Proof of Theorem 1.18.} Suppose that $p\e 1,3,4,5,9\mod {11}$ and so $4p=x^2+11y^2(x,y\in\Bbb Z)$. We choose the sign of $y$ so that $x+y\sqrt{-11}\not\in p\Bbb Z_p$. Taking $D=11$ in [2, Theorem 1.27] and $d_K=-11$ in [7, (12.20)] gives
$$\sum_{k=0}^{p-1}\f{\b{2k}k\b{3k}k\b{6k}{3k}}{(-32)^{3k}}\e \pm \Ls{x+y\sqrt{-11}}2^2\mod {p^3}.$$
 Applying Lemma 2.3 (with $c=4$) we get
$$\sum_{k=0}^{p-1}\f{\b{2k}k\b{3k}k\b{6k}{3k}}{(-32)^{3k}}\e
\pm \f 14\Big(4x^2-2\cdot 4p-\f{4^2p^2}{4x^2}\Big)=\pm\Big(x^2-2p-\f{p^2}{x^2}\Big)\mod{p^3}.$$
By [19, Theorem 4.8],
$\sum_{k=0}^{p-1}\f{\b{2k}k\b{3k}k\b{6k}{3k}}{(-32)^{3k}}
\e\ls {-2}p4x^2\mod p.$
Hence,
$$\sum_{k=0}^{p-1}\f{\b{2k}k\b{3k}k\b{6k}{3k}}{(-32)^{3k}}
\e\Ls {-2}p\Big(x^2-2p-\f{p^2}{x^2}\Big)\mod{p^3}.$$
\par
\par In a similar way, using [2, Theorem 1.27], [7, (12.20)] and [19, Theorems 4.6-4.7 and 4.9] one can prove Theorems 1.16-1.17 and 1.19-1.22.

\section* {\bf 6. Proofs of congruences involving $V_n$}

\par For $\tau\in H$ and $q=e^{2\pi i\tau}$ define
$$s(\tau)=q\prod_{n=1}^{\infty}(1+q^n)^8(1+q^{2n})^8.\eqno{(6.1)}$$
Then clearly
$$s(\tau)=\f{\eta(4\tau)^8}{\eta(\tau)^8}
=\f{f_2(2\tau)^8}{16f_1(2\tau)^8}=\f{f_1(-\f 1{2\tau})^8}{16f_2(-\f 1{2\tau})^8}=\f 1{256s(-\f 1{4\tau})}.\eqno{(6.2)}$$

\par{\bf Lemma 6.1.}
  {\it We have
  $$s\Ls i2=\f 1{16}\q\t{and}\q s\Ls{-1+i}4=-\f 18.$$}
\par{\it Proof.} Since $f_1(-\f 1{\tau})=f_2(\tau)$ we see that
$f_1(i)=f_1(-\f 1i)=f_2(i)$ and so
 $s(\f i2)=\f{f_2(i)^8}{16f_1(i)^8}=\f 1{16}$ by (6.2).
 By the proof of Lemma 4.2,
 $f_2\sls{-1+i}2^8=-4\omega$ and
$f_1\sls{-1+i}2^8=2\omega$.
Thus,
$$s\Ls{-1+i}4=\f{f_2\ls{-1+i}2^8}{16f_1\ls{-1+i}2^8}=\f{-4\omega}{16\cdot 2\omega}=-\f 18,$$
which completes the proof.
\v2
\par{\bf Lemma 6.2}
  {\it The function $s(\tau)$ is a Hauptmodul for $\Gamma_{0}(4)$ with a unique pole at the cusp $[0]$.}
\v2
 \par {\it Proof.}
  For any $\begin{pmatrix}
     a&b\\c&d
 \end{pmatrix}\in \Gamma_{0}(4)$,
 \begin{align*}
     s\left(\frac{a\tau+b}{c\tau+d}\right)
     =\f{\eta(4\frac{a\tau+b}{c\tau+d})^8}{\eta(\frac{a\tau+b}{c\tau+d})^8}
     =\f{\eta(\frac{a(4\tau)+4b}{(c/4)(4\tau)+d})^8}
     {\eta(\frac{a\tau+b}{c\tau+d})^8}.
 \end{align*}
 By (2.1), one finds that
 \begin{align*}
    \eta\Big(\frac{a\tau+b}{c\tau+d}\Big)^8&= \omega^{ab+cd(1-a^{2})-ca+3c_{0}(a-1)+\frac{3}{2}(a^{2}-1)r}
(c\tau+d)^{4}\eta(\tau)^{8},\\
\eta\Big(\frac{a(4\tau)+4b}{(c/4)(4\tau)+d}\Big)^8
&=\omega^{a(4b)+(c/4)d(1-a^{2})-ca+3c_{0}(a-1)+\frac{3}{2}(a^{2}-1)(r-2)}
((c/4)(4\tau)+d)^{4}\eta(4\tau)^{8}
 \end{align*}
  and thus deduces from the preceding line that
  \begin{align*}
      s\Big(\frac{a\tau+b}{c\tau+d}\Big)
     =\f{\eta(\frac{a(4\tau)+4b}{(c/4)(4\tau)+d})^8}
     {\eta(\frac{a\tau+b}{c\tau+d})^8}
     =\omega^{3ab-\frac{9}{2}(a^{2}-1)}\f{\eta(4\tau)^8}{\eta(\tau)^8}
          =\f{\eta(4\tau)^8}{\eta(\tau)^8}=s(\tau).
  \end{align*}
Therefore,  $s(\tau)$ is invariant under $\Gamma_{0}(4)$. Also, it is clear that $s(\tau)$ has no zeros or poles in $H$, and $s(\tau)=q+O(q^{2})$ has a simple zero at the cusp $[i\infty]$, while at the cusp $[0]$ one can find that $s(\tau)$ has a simple pole by noticing via the transformation formula (6.2) that
  $$
  s\Big(-\frac{1}{\tau}\Big)=\f 1{256s\sls {\tau}4}
  =\frac{1}{4^{4}}\cdot\frac{\eta(\tau/4)^{8}}{\eta(\tau)^{8}}
  =\frac{1}{4^{4}}q^{-\frac{1}{4}}+O(q^{-\frac{1}{2}}),
  $$
  where $q=e^{2\pi i\tau}$.
Therefore, by the Residue Theorem for compact Riemann surfaces one can affirm that $s(\tau)$ is a Hauptmodul for $\Gamma_{0}(4)$ with a unique pole at the cusp $[0]$.
\v2
\par{\bf Lemma 6.3}
    {\it  For $\tau\in H$ such that $s(\tau)$ is near~$0$, $s(\tau)$ satisfies that
    \begin{align*}
        \sum_{n=0}^{\infty}V_{n}(-s(\tau))^{n}
        =\frac{\eta(\tau)^{8}}{\eta(2\tau)^{4}}.
    \end{align*}
The eta product on the right hand side is a weight~2 modular form for $\Gamma_{0}(4)$ with a unique zero at the cusp~$[0]$}.
\v2\par
{\it Proof}. Proof of the first claim can be found in [4, Theorem 3.5].
 Next, for any $\begin{pmatrix}
    a&b\\c&d
\end{pmatrix}\in\Gamma_{0}(4)$, by the transformation formula~(2.1), one has that
\begin{align*}
    {\eta\Big(\frac{a\tau+b}{c\tau+d}\Big)^{8}} &= \Big(\frac{a}{c_0}\Big)\zeta_{3}^{ab+cd(1-a^{2})-ca+3c_{0}(a-1)
    +\frac{3}{2}(a^{2}-1)r}
(c\tau+d)^{4}\eta(\tau)^{8},\\
{\eta\Big(2\frac{a\tau+b}{c\tau+d}\Big)^{4}}
&=\eta\Big(\frac{a(2\tau)+2b}{(c/2)(2\tau)+d}\Big)^{4} \\&=\Big(\frac{a}{c_0}\Big)\zeta_{6}^{2ab+\frac{1}{2}cd(1-a^{2})
-\frac{1}{2}ca+3c_{0}(a-1)+\frac{3}{2}(a^{2}-1)(r-1)}
(c\tau+d)^{2}\eta(2\tau)^{4},
\end{align*}
and thus,
$$
 \frac{\eta\left(\frac{a\tau+b}{c\tau+d}\right)^{8}}
 {\eta\left(2\frac{a\tau+b}{c\tau+d}\right)^{4}}
 =(c\tau+d)^{2}\frac{\eta(\tau)^{8}}{\eta(2\tau)^{4}},
$$
that says, $\frac{\eta(\tau)^{8}}{\eta(2\tau)^{4}}$ is a weight~2 modular for~$\Gamma_{0}(4)$.
For the remaining claim, first notice that $\frac{\eta(\tau)^{8}}{\eta(2\tau)^{4}}=1+O(q)$ has no zeros or poles at the cusp $[i\infty]$, and by the transformation formula (2.1) that
$$
\frac{1}{\tau^{2}}\cdot\frac{\eta(-1/\tau)^{8}}{\eta(-2/\tau)^{4}}=-4\frac{\eta(\tau)^{8}}{\eta(\tau/2)^{4}}=-4q^{\frac{1}{4}}+O(q^{\frac{1}{2}}),
$$
which implies that $\frac{\eta(\tau)^{8}}{\eta(2\tau)^{4}}$ has a simple zero at the cusp~$[0]$ of order~1. So one can deduce via the Riemann-Roch theorem that the cusp $[0]$ is the unique zero of $\frac{\eta(\tau)^{8}}{\eta(2\tau)^{4}}$ in $X_{0}(4)$. This completes the proof.

\v2
\par{\bf Proof of Theorem 1.23}. Since $p=x^2+4y^2=(x+2yi)(x-2yi)$, we may choose the sign of $x$ so that $x+2yi\not\in p\Bbb Z_p$.
Let $s(\tau)$ be given by (6.1).
By Lemmas 6.2 and 6.3, one can see that both
$s(\tau)$ and $\sum_{n=0}^{\infty}V_n(-s(\tau))^{n}
$ satisfy the assumptions in Beukers' theorem. Set
 $$a_n=V_n,\q c=4y,\q d=x,\q \alpha=\f{i}2\q\t{and}\q N=4.$$
 Then clearly
 $(c,d)=1,\ N\mid c,$ $c\alpha\bar{\alpha},\ c(\alpha+\bar{\alpha})\in\Bbb Z$,
$p=(c\alpha+d)(c\bar{\alpha}+d)$ and $c\alpha+d\notin p\Bbb Z_p$.
    By Lemma 6.1, $s\sls{i}2=\f 1{16}$.
 Hence, applying Beukers' theorem(i) and Lemma 2.3 we obtain
 \begin{align*}\sum_{n=0}^{p-1}\f{V_n}{(-16)^n}
 &=\sum_{n=0}^{p-1}V_n\Big(-s\Ls{i}2\Big)^n\e \Big(4y\cdot \f{i}2+x\Big)^2
 \\&=(x+y\sqrt{-4})^2\e 4x^2-2p-\f{p^2}{4x^2}\mod {p^3}.\end{align*}
On the other hand, setting
 $$a_n=V_n,\q c=8y,\q d=x+2y,\q \alpha=\f{-1+i}4\q\t{and}\q N=4$$
 one may check that
 $(c,d)=1,\ N\mid c,$ $c\alpha\bar{\alpha},\ c(\alpha+\bar{\alpha})\in\Bbb Z$,
$p=(c\alpha+d)(c\bar{\alpha}+d)$ and $c\alpha+d\notin p\Bbb Z_p$.
     By Lemma 6.1, $s\sls{-1+i}4=-\f 1{8}$.
 Applying Beukers' theorem(i) and Lemma 2.3 we obtain
 \begin{align*}\sum_{n=0}^{p-1}\f{V_n}{8^n}
 &=\sum_{n=0}^{p-1}V_n\Big(-s\ls{-1+i}4\Big)^n\e \Big(8y\cdot \f{-1+i}4+x+2y\Big)^2
 \\&=(x+y\sqrt{-4})^2\e 4x^2-2p-\f{p^2}{4x^2}\mod {p^3}.\end{align*}

\section*{\bf 7. Proofs of congruences involving $T_n$}

\par\q For $\tau\in H$ define
$$w(\tau)=\f{f_2(4\tau)^8}{f_2(\tau)^8}
=\f {f_1(2\tau)^8f_2(4\tau)^8}{16}=\Ls{\eta(\tau)\eta(8\tau)}{\eta(2\tau)
\eta(4\tau)}^8.\eqno{(7.1)}$$
Since $f_1(2\tau)f_2(\tau)=\sqrt 2$ and $f_1(\tau)=f_2(-\f 1{\tau})$ we see that
$$w(\tau)=\f{f_2(4\tau)^8}{f_2(\tau)^8}=\f{16/f_1(8\tau)^8}{16/f_1(2\tau)^8}
=\f{f_1(2\tau)^8}{f_1(8\tau)^8}=\f{f_2(-\f 1{2\tau})^8}{f_2(-\f 1{8\tau})^8}=w\Big(-\f 1{8\tau}\Big).\eqno{(7.2)}$$
\par{\bf Lemma 7.1} {\it We have
$$w\Ls{2+\sqrt{-2}}4=-\f 14,\q w\Ls{5+\sqrt{-7}}{16}=1,\q
w\Ls{1+\sqrt{-7}}8=\f 1{16}.$$}
\par{\it Proof.} Since $f_2(-\f 1{\tau})=f_1(\tau)$ and $f_1(2\tau)f_2(\tau)=\sqrt 2$, we see that
$f_2\sls{\sqrt{-2}}2=f_1(\sqrt {-2})$
and so $f_1(\sqrt{-2})^2=f_1(\sqrt {-2})f_2\sls{\sqrt{-2}}2=\sqrt 2$.
Thus, $f_1(\sqrt{-2})^8=4$. Note that $f_2(2+\sqrt{-2})=\zeta_{24}f_2(1+\sqrt{-2})=\zeta_{24}^2f_2(\sqrt{-2})$
and $f_1\ls{2+\sqrt{-2}}2=\zeta_{48}^{-1}f(\f{\sqrt{-2}}2)
=\zeta_{48}^{-1}f(\sqrt{-2})$. We then have
\begin{align*}&f_2(2+\sqrt{-2})^8=\zeta_{24}^{16}f_2(\sqrt{-2})^8
=\zeta_6^4f_2(\sqrt{-2})^8,
\\&f_1\Ls{2+\sqrt{-2}}2^8=\zeta_{48}^{-8}f(\sqrt{-2})^8=\zeta_6^{-1}
f(\sqrt{-2})^8.
\end{align*}
Hence,
\begin{align*}w\Ls{2+\sqrt{-2}}4&=\f{f_1\sls{2+\sqrt{-2}}2^8f_2(2+\sqrt {-2})^8}{16}
=\f{\zeta_6^3f(\sqrt{-2})^8f_2(\sqrt{-2})^8}{16}
\\&=-\f{(\sqrt 2)^8}{16f_1(\sqrt{-2})^8}=-\f 14.\end{align*}
\par Using the properties of $f(\tau),f_1(\tau)$ and $f_2(\tau)$ we derive that
$$f_1\Ls{5+\sqrt{-7}}2=\zeta_{48}^{-1}f\Ls{3+\sqrt{-7}}2
=\zeta_{48}^{-2}f_1\Ls{1+\sqrt{-7}}2$$
and
\begin{align*}f_1\Ls{5+\sqrt{-7}}8
&=\zeta_{48}^{-1}f\Ls{-3+\sqrt{-7}}8=
\zeta_{48}^{-1}f\Big(-\f 8{-3+\sqrt{-7}}\Big)
\\&=\zeta_{48}^{-1}f\Ls{3+\sqrt{-7}}2
=\zeta_{48}^{-2}f_1\Ls{1+\sqrt{-7}}2.
\end{align*}
Thus, $f_1\sls{5+\sqrt{-7}}2=f_1\sls{5+\sqrt{-7}}8$ and so
$$w\Ls{5+\sqrt{-7}}{16}=\f{f_1\sls{5+\sqrt{-7}}8^8}
{f_1\sls{5+\sqrt{-7}}2^8}=1.$$
\par By the proof of Lemma 4.2, $f_2\sls{1+\sqrt{-7}}2^8=-\omega^2$. Thus, applying Lemma 2.1 gives $$f_2\Ls{-1+\sqrt{-7}}2^8=\overline{f_2\Ls{1+\sqrt{-7}}2}^8
=\overline{f_2\Ls{1+\sqrt{-7}}2^8}=\overline{-\omega^2}=-\omega.$$
Since $\f{1+\sqrt{-7}}4\cdot\f{-1+\sqrt{-7}}2=-1$, we have
$f_1\sls{1+\sqrt{-7}}4^8=f_2\sls{-1+\sqrt{-7}}2^8=-\omega$. Therefore,
$$w\Ls{1+\sqrt{-7}}8=\f{f_1\sls{1+\sqrt{-7}}4^8
f_2\sls{1+\sqrt{-7}}2^8}{16}=\f{-\omega\cdot(-\omega^2)}{16}=\f 1{16}.$$
\par{\bf Remark 7.1} Lemma 7.1 was stated by Beukers in [2, p.32] without proof.
\v2
\par{\bf Lemma 7.2} {\it The function $w(\tau)$ is a Hauptmodul for $\Gamma_{0}(8)^{+}$ with a unique pole at the cusp~$[\frac{1}{2}]$ and satisfies that
$$
\sum_{n=0}^{\infty}T_{n}w(\tau)^{n}=\frac{\eta(2\tau)^{6}\eta(4\tau)^{6}}{\eta(\tau)^{4}\eta(8\tau)^{4}},
$$
and the eta quotient $\frac{\eta(2\tau)^{6}\eta(4\tau)^{6}}{\eta(\tau)^{4}\eta(8\tau)^{4}}$ is a weight~2 modular form with respect to~$\Gamma_{0}(8)^{+}$ with a unique zero at the cusp~$[\frac{1}{2}]$.
}\v2
\par{\it Proof.}
 For any $\begin{pmatrix}
    a&b\\c&d
\end{pmatrix}\in\Gamma_{0}(8)$, using the transformation formula~(2.1), one finds that
\begin{align*}
 {\eta\left(\frac{a\tau+b}{c\tau+d}\right)^{8}} &= \omega^{ab+cd(1-a^{2})-ca+3c_{0}(a-1)+r\frac{3}{2}(a^{2}-1)}
(c\tau+d)^{4}\eta(\tau)^{8},\\
 {\eta\left(2\frac{a\tau+b}{c\tau+d}\right)^{8}} &=\eta\left(\frac{a(2\tau)+2b}{(c/2)(2\tau)+d}\right)^{8}\\ &=\omega^{2ab+\frac{1}{2}cd(1-a^{2})-ca+3c_{0}(a-1)+(r-1)\frac{3}{2}(a^{2}-1)}
(c\tau+d)^{4}\eta(2\tau)^{8},\\
 {\eta\left(4\frac{a\tau+b}{c\tau+d}\right)^{8}} &=\eta\left(\frac{a(4\tau)+4b}{(c/4)(4\tau)+d}\right)^{8}\\ &=\omega^{4ab+\frac{1}{4}cd(1-a^{2})-ca+3c_{0}(a-1)+(r-2)\frac{3}{2}(a^{2}-1)}
(c\tau+d)^{4}\eta(4\tau)^{8},\\
 {\eta\left(8\frac{a\tau+b}{c\tau+d}\right)^{8}} &=\eta\left(\frac{a(8\tau)+8b}{(c/8)(8\tau)+d}\right)^{8}\\ &=\omega^{8ab+\frac{1}{8}cd(1-a^{2})-ca+3c_{0}(a-1)+(r-3)\frac{3}{2}(a^{2}-1)}
(c\tau+d)^{4}\eta(8\tau)^{8},\\
\end{align*}
and thus,
$$
w\left(\frac{a\tau+b}{c\tau+d}\right)=w(\tau).
$$
By (7.2), $w(-\frac{1}{8\tau})=w(\tau)$.
Therefore, $w(\tau)$ is invariant under $\Gamma_{0}(8)^{+}$. Since $w(\tau)$ is an infinite product in $H$ and $w(\tau)=q+O(q^{2})$, where $q=e^{2\pi i\tau}$, $w(\tau)$ has a simple zero at the cusp~$[i\infty]$, and so has a unique pole at the other cusp~$[\frac{1}{2}]$. Similarly,
\begin{align*}
 {\eta\left(\frac{a\tau+b}{c\tau+d}\right)^{4}} &= \zeta_{6}^{ab+cd(1-a^{2})-ca+3c_{0}(a-1)+r\frac{3}{2}(a^{2}-1)}
(c\tau+d)^{2}\eta(\tau)^{4},\\
 {\eta\left(2\frac{a\tau+b}{c\tau+d}\right)^{6}} &=\eta\left(\frac{a(2\tau)+2b}{(c/2)(2\tau)+d}\right)^{6}\\ &=i^{2ab+\frac{1}{2}cd(1-a^{2})-ca+3c_{0}(a-1)+(r-1)\frac{3}{2}(a^{2}-1)}
(c\tau+d)^{3}\eta(2\tau)^{6},
\end{align*}
\begin{align*}
 {\eta\left(4\frac{a\tau+b}{c\tau+d}\right)^{6}} &=\eta\left(\frac{a(4\tau)+4b}{(c/4)(4\tau)+d}\right)^{6}\\ &=i^{4ab+\frac{1}{4}cd(1-a^{2})-ca+3c_{0}(a-1)+(r-2)\frac{3}{2}(a^{2}-1)}
(c\tau+d)^{3}\eta(4\tau)^{6},\\
 {\eta\left(8\frac{a\tau+b}{c\tau+d}\right)^{4}} &=\eta\left(\frac{a(8\tau)+8b}{(c/8)(8\tau)+d}\right)^{4}\\ &=\zeta_{6}^{8ab+\frac{1}{8}cd(1-a^{2})-ca+3c_{0}(a-1)+(r-3)\frac{3}{2}(a^{2}-1)}
(c\tau+d)^{2}\eta(8\tau)^{4},\\
\end{align*}
and thus after some simple cancellations, one finds that
$$
\frac{\eta(2\frac{a\tau+b}{c\tau+d})^{6}\eta(4\frac{a\tau+b}{c\tau+d})^{6}}{\eta(\frac{a\tau+b}{c\tau+d})^{4}\eta(8\frac{a\tau+b}{c\tau+d})^{4}}=(c\tau+d)^{2}\frac{\eta(2\tau)^{6}\eta(4\tau)^{6}}{\eta(\tau)^{4}\eta(8\tau)^{4}}.
$$
Also, under $\tau\to-1/8\tau$, one can find that
$$
\frac{\eta(-1/4\tau)^{6}\eta(1/2\tau)^{6}}{\eta(1/8\tau)^{4}\eta(1/\tau)^{4}}=\frac{(4\tau/i)^{3}\eta(4\tau)^{6}(2\tau/i)^{3}\eta(2\tau)^{6}}{(8\tau/i)^{2}\eta(\tau)^{4}(\tau/i)^{2}\eta(8\tau)^{4}}=-\frac{\eta(2\tau)^{6}\eta(4\tau)^{6}}{\eta(\tau)^{4}\eta(8\tau)^{4}},
$$
which together with the preceding result implies that
$\frac{\eta(2\tau)^{6}\eta(4\tau)^{6}}{\eta(\tau)^{4}\eta(8\tau)^{4}}$ is a weight~2 modular form for $\Gamma_{0}(8)^{+}$. Clearly, by its infinite product representation it has no zeros or poles at any point of $X_{0}(8)^{+}$ apart from the cusp $[\frac{1}{2}]$, so by the Riemann-Roch theorem, the eta quotient $\frac{\eta(2\tau)^{6}\eta(4\tau)^{6}}{\eta(\tau)^{4}\eta(8\tau)^{4}}$ must have a unique zero at the cusp~$[\frac{1}{2}]$. Proof of the relation between $w(\tau)$ and the eta quotient $\frac{\eta(2\tau)^{6}\eta(4\tau)^{6}}{\eta(\tau)^{4}\eta(8\tau)^{4}}$ can be found in [5, Section~8.8].
 \v2
\par{\bf Proof of Theorem 1.24.} Since $p=x^2+2y^2=(x+y\sqrt{-2})(x-y\sqrt{-2})$, we may choose the sign of $x$ so that $x+y\sqrt{-2}\not\in p\Bbb Z_p$. Let $w(\tau)$ be given by (7.1). By Lemma 7.2, one can see that both
$w(\tau)$ and $\sum_{n=0}^{\infty}T_nw(\tau)^{n}
$ satisfy the assumptions in Beukers' theorem.
\par For $p=x^2+2y^2\e 1\mod 8$ we have $2\mid y$. Set
 $$a_n=T_n,\q c=4y,\q d=x-2y,\q \alpha=\f{2+\sqrt{-2}}4\q\t{and}\q N=8.$$
 Then clearly
 $(c,d)=1,\ N\mid c,$ $c\alpha\bar{\alpha}, c(\alpha+\bar{\alpha})\in\Bbb Z$,
$p=(c\alpha+d)(c\bar{\alpha}+d)$ and $c\alpha+d\notin p\Bbb Z_p$.
    By Lemma 7.1, $w\sls{2+\sqrt{-2}}4=-\f 1{4}$.
 Hence, applying Beukers' theorem(i) and Lemma 2.3 we obtain
 \begin{align*}\sum_{n=0}^{p-1}\f{T_n}{(-4)^n}
 &=\sum_{n=0}^{p-1}T_nw\Ls{2+\sqrt{-2}}4^n\e \Big(4y\cdot \f{2+\sqrt{-2}}4+x-2y\Big)^2
 \\&=(x+y\sqrt{-2})^2\e 4x^2-2p-\f{p^2}{4x^2}\mod {p^3}.\end{align*}

\par For $p=x^2+2y^2\e 3\mod 8$ we have $2\nmid xy$.
Set $$a_n=T_n,\q
c=x,\q d=-\f{x+y}2,\q \alpha=\f{2+\sqrt{-2}}4\q\t{and}\q N=8.$$
It is easy to check that $(c,dN)=1,\ Nc\alpha\bar{\alpha}\in\Bbb Z,\ $  $c(\alpha+\bar{\alpha})\in\Bbb Z$,
$p=N(c\alpha+d)(c\bar{\alpha}+d)$ and $c\alpha+d\notin p\Bbb Z_p$.
Since $w\ls{2+\sqrt{-2}}4=-\f 1{4}$ by Lemma 7.1, using Beukers' theorem(ii) and Lemma 2.3 we see that
 \begin{align*}\sum_{n=0}^{p-1}\f{T_n}{(-4)^n}
 &=\sum_{n=0}^{p-1}T_n w\Ls{2+\sqrt{-2}}4^n\e
 -8\Big(x\cdot\f{2+\sqrt{-2}}4-\f{x+y}2\Big)^2\\&
  =(x+y\sqrt {-2})^2\e 4x^2-2p-\f{p^2}{4x^2}\mod {p^3}.\end{align*}
\v2
\par{\bf Proof of Theorem 1.25.} Since $p=x^2+7y^2=(x+y\sqrt{-7})(x-y\sqrt{-7})$, we may choose the sign of $y$ so that $x+y\sqrt{-7}\not\in p\Bbb Z_p$. Let $w(\tau)$ be given by (7.1). By Lemma 7.2, one can see that both
$w(\tau)$ and $\sum_{n=0}^{\infty}T_nw(\tau)^{n}
$ satisfy the assumptions in Beukers' theorem.
\par Since $x\not\e y\mod 2$, setting
 $$a_n=T_n,\q c=16y,\q d=x-5y,\q \alpha=\f{5+\sqrt{-7}}{16}\q\t{and}\q N=8,$$
 we see that
 $(c,d)=1,\ N\mid c,\;$ $c\alpha\bar{\alpha}, c(\alpha+\bar{\alpha})\in\Bbb Z$,
$p=(c\alpha+d)(c\bar{\alpha}+d)$ and $c\alpha+d\notin p\Bbb Z_p$.
    By Lemma 7.1, $w\sls{5+\sqrt{-7}}{16}=1$.
 Hence, applying Beukers' theorem(i) and Lemma 2.3 we obtain
 \begin{align*}\sum_{n=0}^{p-1}T_n
 &=\sum_{n=0}^{p-1}T_nw\Ls{5+\sqrt{-7}}{16}^n\e \Big(16y\cdot \f{5+\sqrt{-7}}{16}+x-5y\Big)^2
 \\&=(x+y\sqrt{-7})^2\e 4x^2-2p-\f{p^2}{4x^2}\mod {p^3}.\end{align*}
  On the other hand, setting
 $$a_n=T_n,\q c=8y,\q d=x-y,\q
 \alpha=\f{1+\sqrt{-7}}{8}\q\t{and}\q N=8,$$
 we find that
 $(c,d)=1,\ N\mid c,\;$ $c\alpha\bar{\alpha}, c(\alpha+\bar{\alpha})\in\Bbb Z$,
$p=(c\alpha+d)(c\bar{\alpha}+d)$ and $c\alpha+d\notin p\Bbb Z_p$.
    By Lemma 7.1, $w\sls{1+\sqrt{-7}}{8}=\f 1{16}$.
 Hence, applying Beukers' theorem(i) and Lemma 2.3 yields \begin{align*}\sum_{n=0}^{p-1}\f{T_n}{16^n}
 &=\sum_{n=0}^{p-1}T_nw\Ls{1+\sqrt{-7}}{8}^n\e \Big(8y\cdot \f{1+\sqrt{-7}}{8}+x-y\Big)^2
 \\&=(x+y\sqrt{-7})^2\e 4x^2-2p-\f{p^2}{4x^2}\mod {p^3}.\end{align*}
\v2

\section*{8. Proofs of congruences involving $D_n$}
\par For $\tau\in H$ let
$$v(\tau)=\f 1{f(2\tau)^6f(6\tau)^6}.\eqno{(8.1)}$$
Since $f(\tau)=f(-\f 1{\tau})$ we see that
$$v(\tau)=\f 1{f(-\f 1{2\tau})^6f(-\f 1{6\tau})^6}=v\Big(-\f 1{12\tau}\Big).\eqno{(8.2)}$$
As $f(\tau)f_1(\tau)f_2(\tau)=\sqrt 2$, we also have
$$v(\tau)=\f 1{2^6}\big(f_1(2\tau)f_2(2\tau)f_1(6\tau)f_2(6\tau)\big)^6
=\Ls{\eta(\tau)\eta(3\tau)\eta(4\tau)\eta(12\tau)}
{\eta(2\tau)^2\eta(6\tau)^2}^6.\eqno{(8.3)}$$
\par{\bf Lemma 8.1} {\it We have
\begin{align*}& v\Ls{1+\sqrt{-2}}6=\f 18,\q v\Ls{3+\sqrt{-3}}6=-\f 12\q v\Ls{3+2\sqrt{-3}}6=-\f 1{32},
\\& v\Ls{3+\sqrt{-6}}6=-\f 18.\end{align*}}
\par{\it Proof.} Since $f(\tau+1)=\zeta_{48}^{-1}f_1(\tau)$, $f_1(\tau+1)=\zeta_{48}^{-1}f(\tau)$ and $f(\tau)=f(-\f 1{\tau})$, we see that
$f(1+\sqrt{-2})=\zeta_{48}^{-1}f_1(\sqrt{-2})$ and
$f\sls{1+\sqrt{-2}}3=f(-1+\sqrt{-2})=\zeta_{48}f_1(\sqrt{-2})$. Thus,
$f(1+\sqrt{-2})f\sls{1+\sqrt{-2}}3=f_1(\sqrt{-2})^2$. By the proof of Lemma 3.2, $f_1(\sqrt{-2})^2=\sqrt{2}$. Hence,
$$v\Ls{1+\sqrt{-2}}6=\f 1{f(1+\sqrt{-2})^6f\sls{1+\sqrt{-2}}3^6}
=\f 1{f_1(\sqrt{-2})^{12}}=\f 1{(\sqrt 2)^6}=\f 18.$$
\par From the properties of Weber's functions, we see that
\begin{align*}&f(3+\sqrt{-3})=\zeta_{48}^{-1}f_1(2+\sqrt{-3})=\zeta_{48}^{-2}f(1+\sqrt{-3})
=\zeta_{48}^{-3}f_1(\sqrt{-3}),
\\&f\Ls{3+\sqrt{-3}}3=f\Big(1-\f 1{\sqrt{-3}}\Big)=\zeta_{48}^{-1}
f_1\Big(-\f 1{\sqrt{-3}}\Big)=\zeta_{48}^{-1}f_2(\sqrt{-3}).
\end{align*}
Thus,
\begin{align*}v\Ls{3+\sqrt{-3}}6&=\f 1{f\sls{3+\sqrt{-3}}3^6f(3+\sqrt{-3})^6}
=\f 1{\zeta_{48}^{-6}f_2(\sqrt{-3})^6\cdot\zeta_{48}^{-18}f_1(\sqrt{-3})^6}
\\&=\zeta_{48}^{24}\Ls{f(\sqrt{-3})}{\sqrt 2}^6=-\f {f(\sqrt{-3})^6}8.\end{align*}
By Lemma 3.2, $f(\sqrt{-3})^{24}=t(\f{\sqrt{-3}}2)^{-1}=256=4^4$. Thus,
$f(\sqrt{-3})^6\in\{4,-4,4i,-4i\}$. Since
$$f(\sqrt{-3})^6=e^{-2\pi i\sqrt{-3}/8}\prod_{n=1}^{\infty}
\big(1+e^{2\pi i\sqrt{-3}(n-\f 12)}\big)^6
=e^{2\pi \sqrt{3}/8}\prod_{n=1}^{\infty}
\big(1+e^{-2\pi \sqrt{3}(n-\f 12)}\big)^6,$$
we see that $f(\sqrt{-3})^6$ is real and positive. Hence $f(\sqrt{-3})^6=4$ and so
$$v\Ls{3+\sqrt{-3}}6=-\f {f(\sqrt{-3})^6}8=-\f 48=-\f 12.$$
\par Using the properties of Weber functions, we find that
$$f\ls{3+2\sqrt{-3}}3=f\Big(1-\f 2{\sqrt{-3}}\Big)=\zeta_{48}^{-1}f_1
\Big(-\f 2{\sqrt{-3}}\Big)=\zeta_{48}^{-1}f_2\Ls{\sqrt{-3}}2
=\zeta_{48}^{-1}\f{\sqrt 2}{f_1(\sqrt{-3})}$$
and
$$f(3+2\sqrt{-3})=\zeta_{48}^{-1}f_1(2+2\sqrt{-3})=\zeta_{48}^{-2}
f(1+2\sqrt{-3})=\zeta_{48}^{-3}f_1(2\sqrt{-3})=\zeta_{48}^{-3}\f{\sqrt 2}{f_2(\sqrt {-3})}.$$
Thus,
$$f\Ls{3+2\sqrt{-3}}3f(3+2\sqrt{-3})=\zeta_{48}^{-4}\cdot\f{\sqrt 2\cdot\sqrt 2}{f_1(\sqrt{-3})f_2(\sqrt{-3})}=\zeta_{12}^{-1}\cdot\sqrt 2f(\sqrt {-3})$$
and therefore
$$v\Ls{3+2\sqrt{-3}}6=\f 1{f\sls{3+2\sqrt{-3}}3^6f(3+2\sqrt{-3})^6}
=\f 1{\zeta_{12}^{-6}(\sqrt 2)^6f(\sqrt{-3})^6}=-\f 1{8f(\sqrt{-3})^6}.$$
By the previous argument, $f(\sqrt{-3})^6=4$. Hence $v\sls{3+2\sqrt{-3}}6=-\f 1{8f(\sqrt{-3})^6}=-\f 1{32}.$
\par Using the properties of Weber's functions, we see that
$f\ls{3+\sqrt{-6}}3=\zeta_{48}^{-1}f_1\ls{\sqrt{-6}}3=\zeta_{48}^{-1}
f_2\ls{\sqrt{-6}}2$ and $f(3+\sqrt{-6})=\zeta_{48}^{-1}f_1(2+\sqrt{-6})
=\zeta_{48}^{-2}f(1+\sqrt{-6})=\zeta_{48}^{-3}f_1(\sqrt{-6})$. Thus,
\begin{align*}v\Ls{3+\sqrt{-6}}6&=\f 1{f\sls{3+\sqrt{-6}}3^6f(3+\sqrt{-6})^6}
=\f 1{(\zeta_{48}^{-1}
f_2\ls{\sqrt{-6}}2\cdot \zeta_{48}^{-3}f_1(\sqrt{-6}))^6}
\\&=\f 1{(\zeta_{48}^{-4}\sqrt 2)^6}=-\f 18.\end{align*}

\v2
\par{\bf Remark 8.1} Lemma 8.1 was stated by Beukers in [2, p.31] without proof.
\v2
\par{\bf Lemma 8.2} {\it The function $v(\tau)$ is a Hauptmodul for $\Gamma_{0}(12)^{+}$ with a unique pole at the cusp~$[\frac{1}{2}]$ and satisfies that
$$
\sum_{n=0}^{\infty}D_{n}v(\tau)^{n}=\frac{\eta(2\tau)^{10}\eta(6\tau)^{10}}{\eta(\tau)^{4}\eta(3\tau)^{4}\eta(4\tau)^{4}\eta(12\tau)^{4}}
$$
for $\tau$ such that $v(\tau)$ is near~0. In particular, the eta product on the right hand side is a weight~2 modular form for $\Gamma_{0}(12)^{+}$ with a unique zero at the cusp~$[\frac{1}{2}]$.}
\v2
\par {\it Proof.}
For any $\begin{pmatrix}
    a&b\\c&d
\end{pmatrix}\in \Gamma_{0}(12)$, by (2.1), one has that
\begin{align*}
\eta\Big(\frac{a\tau+b}{c\tau+d}\Big)^{6}
&=i^{ab+cd(1-a^{2})-ca+3c_{0}(a-1)+\frac{3}{2}(a^{2}-1)r}
(c\tau+d)^{3}\eta(\tau)^{6},\\
\eta\Big(3\frac{a\tau+b}{c\tau+d}\Big)^{6}
&=\eta\Big(\frac{a(3\tau)+3b}{(c/3)(3\tau)+d}\Big)^{6}
\\&=i^{3ab+\frac{1}{3}cd(1-a^{2})-\frac{1}{3}ca+3c_{0}(a-1)
+\frac{3}{2}(a^{2}-1)r}
(c\tau+d)^{3}\eta(3\tau)^{6},\\
\eta\Big(4\frac{a\tau+b}{c\tau+d}\Big)^{6}
&=\eta\Big(\frac{a(4\tau)+4b}{(c/4)(4\tau)+d}\Big)^{6}
\\&=i^{4ab+\frac{1}{4}cd(1-a^{2})-\frac{1}{4}ca+3c_{0}(a-1)
+\frac{3}{2}(a^{2}-1)(r-2)}
(c\tau+d)^{3}\eta(4\tau)^{6},\\
\eta\Big(12\frac{a\tau+b}{c\tau+d}\Big)^{6}
&=\eta\Big(\frac{a(12\tau)+12b}{(c/12)(12\tau)+d}\Big)^{6}
\\&=i^{12ab+\frac{1}{12}cd(1-a^{2})-\frac{1}{12}ca+3c_{0}(a-1)
+\frac{3}{2}(a^{2}-1)(r-2)}
(c\tau+d)^{3}\eta(12\tau)^{6},\\
\eta\Big(2\frac{a\tau+b}{c\tau+d}\Big)^{12}
&=\eta\Big(\frac{a(2\tau)+2b}{(c/2)(2\tau)+d}\Big)^{12}
\\&=(-1)^{2ab+\frac{1}{2}cd(1-a^{2})-\frac{1}{2}ca+3c_{0}(a-1)
+\frac{3}{2}(a^{2}-1)(r-1)}
(c\tau+d)^{6}\eta(2\tau)^{12},\\
\eta\Big(6\frac{a\tau+b}{c\tau+d}\Big)^{12}
&=\eta\Big(\frac{a(6\tau)+6b}{(c/6)(6\tau)+d}\Big)^{12}
\\&=(-1)^{6ab+\frac{1}{6}cd(1-a^{2})-\frac{1}{6}ca+3c_{0}(a-1)
+\frac{3}{2}(a^{2}-1)(r-1)}
(c\tau+d)^{6}\eta(6\tau)^{12}.
\end{align*}
Putting these together and simplifying the resulting multiplier, one shows that
$$
v\left(\frac{a\tau+b}{c\tau+d}\right)=v(\tau)
$$
for any $\begin{pmatrix}
    a&b\\c&d
\end{pmatrix}\in\Gamma_{0}(12)$. Also,
$v(-\f 1{12\tau})=v(\tau)$ by (8.2).
Therefore, $v(\tau)$ is invariant under $\Gamma_{0}(12)^{+}$ with a unique zero at the cusp~$[i\infty]$ and a unique pole at the cusp~$[\frac{1}{2}]$, and so it is a Hauptmodul for $\Gamma_{0}(12)^{+}$.  Proof of the eta product representation for $\sum_{n=0}^{\infty}D_{n}v(\tau)^{n}$ can be found in [5, Section 12.4]. Also, for $\begin{pmatrix}
    a&b\\c&d
\end{pmatrix}\in\Gamma_{0}(12)$,
\begin{align*}
\eta\left(\frac{a\tau+b}{c\tau+d}\right)^{4}
&=\zeta_{6}^{ab+cd(1-a^{2})-ca+3c_{0}(a-1)+r\frac{3}{2}(a^{2}-1)}
(c\tau+d)^{2}\eta(\tau)^{4},\\
\eta\left(3\frac{a\tau+b}{c\tau+d}\right)^{4}
&=\eta\left(\frac{a(3\tau)+3b}{(c/3)(3\tau)+d}\right)^{4}
\\&=\zeta_{6}^{3ab+\frac{1}{3}cd(1-a^{2})-\frac{1}{3}ca+3c_{0}(a-1)+r\frac{3}{2}(a^{2}-1)}
(c\tau+d)^{2}\eta(3\tau)^{4},\\
\eta\left(4\frac{a\tau+b}{c\tau+d}\right)^{4}
&=\eta\left(\frac{a(4\tau)+4b}{(c/4)(4\tau)+d}\right)^{4}
\\&=\zeta_{6}^{4ab+\frac{1}{4}cd(1-a^{2})-\frac{1}{4}ca+3c_{0}(a-1)+(r-2)\frac{3}{2}(a^{2}-1)}
(c\tau+d)^{2}\eta(4\tau)^{4},\\
\eta\left(12\frac{a\tau+b}{c\tau+d}\right)^{4}
&=\eta\left(\frac{a(12\tau)+12b}{(c/12)(12\tau)+d}\right)^{4}
\\&=\zeta_{6}^{12ab+\frac{1}{12}cd(1-a^{2})-\frac{1}{12}ca+3c_{0}(a-1)+(r-2)\frac{3}{2}(a^{2}-1)}
(c\tau+d)^{2}\eta(12\tau)^{4},\\
\eta\left(2\frac{a\tau+b}{c\tau+d}\right)^{10}
&=\eta\left(\frac{a(2\tau)+2b}{(c/2)(2\tau)+d}\right)^{10}
\\&=\zeta_{12}^{5(2ab+\frac{1}{2}cd(1-a^{2})-\frac{1}{2}ca+3c_{0}(a-1)+(r-1)\frac{3}{2}(a^{2}-1))}
(c\tau+d)^{5}\eta(2\tau)^{10},\\
\eta\left(6\frac{a\tau+b}{c\tau+d}\right)^{10}
&=\eta\left(\frac{a(6\tau)+6b}{(c/6)(6\tau)+d}\right)^{10}
\\&=\zeta_{12}^{5(6ab+\frac{1}{6}cd(1-a^{2})-\frac{1}{6}ca+3c_{0}(a-1)+(r-1)\frac{3}{2}(a^{2}-1))}
(c\tau+d)^{5}\eta(6\tau)^{10},
\end{align*}
so that the eta product $F(\tau)=\frac{\eta(2\tau)^{10}\eta(6\tau)^{10}}{\eta(\tau)^{4}\eta(3\tau)^{4}\eta(4\tau)^{4}\eta(12\tau)^{4}}$
actually satisfies that
$$
F\left(\frac{a\tau+d}{c\tau+d}\right)=(c\tau+d)^{2}F(\tau).$$
Similarly, one finds that
$$\frac{\eta(1/6\tau)^{10}\eta(1/2\tau)^{10}}{\eta(-1/12\tau)^{4}\eta(-1/4\tau)^{4}\eta(-1/3\tau)^{4}\eta(-1/\tau)^{4}}=-\tau^{2}\frac{\eta(2\tau)^{10}\eta(6\tau)^{10}}{\eta(\tau)^{4}\eta(3\tau)^{4}\eta(4\tau)^{4}\eta(12\tau)^{4}}.$$ Therefore, the eta product $\frac{\eta(2\tau)^{10}\eta(6\tau)^{10}}{\eta(\tau)^{4}\eta(3\tau)^{4}\eta(4\tau)^{4}\eta(12\tau)^{4}}$ is a weight~2 modular form for~$\Gamma_{0}(12)^{+}$ with a unique zero at the cusp~$[\frac{1}{2}]$.

\v2
\par Let $v(\tau)$ be given by (8.1). By Lemma 8.2, one can see that both $v(\tau)$ and $\sum_{n=0}^{\infty}D_nv(\tau)^{n}$ satisfy the assumptions in Beukers' theorem.
\v2
\par{\bf Proof of Theorem 1.26.} Since $p=x^2+2y^2=(x+y\sqrt{-2})(x-y\sqrt{-2})$, we may choose the sign of $x$ so that $x+y\sqrt{-2}\not\in p\Bbb Z_p$.
Since $p=x^2+2y^2\e 1\mod 8$ we have $2\mid y$. Set
 $$a_n=D_n,\q c=6y,\q d=x-y,\q \alpha=\f{1+\sqrt{-2}}6\q\t{and}\q N=12.$$
 Then clearly
 $(c,d)=1,\ N\mid c,$ $c\alpha\bar{\alpha}, c(\alpha+\bar{\alpha})\in\Bbb Z$,
$p=(c\alpha+d)(c\bar{\alpha}+d)$ and $c\alpha+d\notin p\Bbb Z_p$.
    By Lemma 8.1, $v\sls{1+\sqrt{-2}}6=\f 1{8}$.
 Hence, applying Beukers' theorem(i) and Lemma 2.3 we obtain
 \begin{align*}\sum_{n=0}^{p-1}\f{D_n}{8^n}
 &=\sum_{n=0}^{p-1}D_nv\Ls{1+\sqrt{-2}}6^n\e \Big(6y\cdot \f{1+\sqrt{-2}}6+x-y\Big)^2
 \\&=(x+y\sqrt{-2})^2\e 4x^2-2p-\f{p^2}{4x^2}\mod {p^3}.\end{align*}

\par{\bf Proof of Theorem 1.27.} Since $p=x^2+3y^2=(x+y\sqrt{-3})(x-y\sqrt{-3})$, we may choose the sign of $x$ so that $x+y\sqrt{-3}\not\in p\Bbb Z_p$.
 Since $p=x^2+3y^2\e 1\mod {12}$ we have $2\nmid x$ and $2\mid y$. Set
 $$a_n=D_n,\q c=6y,\q d=x-3y,\q \alpha=\f{3+\sqrt{-3}}6\q\t{and}\q N=12.$$
 Then clearly
 $(c,d)=1,\ N\mid c,$ $c\alpha\bar{\alpha}, c(\alpha+\bar{\alpha})\in\Bbb Z$,
$p=(c\alpha+d)(c\bar{\alpha}+d)$ and $c\alpha+d\notin p\Bbb Z_p$.
    By Lemma 8.1, $v\sls{3+\sqrt{-3}}6=-\f 1{2}$.
 Hence, applying Beukers' theorem(i) and Lemma 2.3 we obtain
 \begin{align*}\sum_{n=0}^{p-1}\f{D_n}{(-2)^n}
 &=\sum_{n=0}^{p-1}D_nv\Ls{3+\sqrt{-3}}6^n\e \Big(6y\cdot \f{3+\sqrt{-3}}6+x-3y\Big)^2
 \\&=(x+y\sqrt{-3})^2\e 4x^2-2p-\f{p^2}{4x^2}\mod {p^3}.\end{align*}
Similarly, for $p=x^2+3y^2\e 1\mod {24}$ we have $2\nmid x$ and $4\mid y$. On setting
 $$a_n=D_n,\q c=3y,\q d=x-\f{3y}2,\q \alpha=\f{3+2\sqrt{-3}}6\q\t{and}\q N=12$$
 we find that
 $(c,d)=1,\ N\mid c,$ $c\alpha\bar{\alpha}, c(\alpha+\bar{\alpha})\in\Bbb Z$,
$p=(c\alpha+d)(c\bar{\alpha}+d)$ and $c\alpha+d\notin p\Bbb Z_p$.
    By Lemma 8.1, $v\sls{3+2\sqrt{-3}}6=-\f 1{32}$.
 Hence, applying Beukers' theorem(i) and Lemma 2.3 we obtain
 \begin{align*}\sum_{n=0}^{p-1}\f{D_n}{(-32)^n}
 &=\sum_{n=0}^{p-1}D_nv\Ls{3+2\sqrt{-3}}6^n\e \Big(3y\cdot \f{3+2\sqrt{-3}}6+x-\f{3y}2\Big)^2
 \\&=(x+y\sqrt{-3})^2\e 4x^2-2p-\f{p^2}{4x^2}\mod {p^3}.\end{align*}

\par{\bf Proof of Theorem 1.28.} For $p=x^2+6y^2\e 1\mod {24}$ we have $2\mid y$. One may choose the sign of $x$ so that $x+y\sqrt{-6}\notin p\Bbb Z_p$. Setting
 $$a_n=D_n,\q c=6y,\q d=x-3y,\q \alpha=\f{3+\sqrt{-6}}6\q\t{and}\q N=12$$
 we find that
 $(c,d)=1,\ N\mid c,$ $c\alpha\bar{\alpha}, c(\alpha+\bar{\alpha})\in\Bbb Z$,
$p=(c\alpha+d)(c\bar{\alpha}+d)$ and $c\alpha+d\notin p\Bbb Z_p$.
    By Lemma 8.1, $v\sls{3+\sqrt{-6}}6=-\f 1{8}$.
 Hence, applying Beukers' theorem(i) and Lemma 2.3 we obtain
 \begin{align*}\sum_{n=0}^{p-1}\f{D_n}{(-8)^n}
 &=\sum_{n=0}^{p-1}D_nv\Ls{3+\sqrt{-6}}6^n\e \Big(6y\cdot \f{3+\sqrt{-6}}6+x-3y\Big)^2
 \\&=(x+y\sqrt{-6})^2\e 4x^2-2p-\f{p^2}{4x^2}\mod {p^3}.
\end{align*}
\par For $p=2x^2+3y^2\e 5\mod {24}$ we have $2\nmid x$ and $2p=(2x)^2+6y^2$. One may choose the sign of $x$ so that $2x+y\sqrt{-6}\notin p\Bbb Z_p$. Setting
 $$a_n=D_n,\q c=x,\q d=-\f{x+y}2,\q \alpha=\f{3+\sqrt{-6}}6\q\t{and}\q N=12$$
 we find that
 $(c,dN)=1,\ Nc\alpha\bar{\alpha}\in\Bbb Z,\ $  $c(\alpha+\bar{\alpha})\in\Bbb Z$,
$p=N(c\alpha+d)(c\bar{\alpha}+d)$ and $c\alpha+d\notin p\Bbb Z_p$.
Since $v\ls{3+\sqrt{-6}}6=-\f 1{8}$ by Lemma 8.1, using Beukers' theorem(ii) and Lemma 2.3 we see that
 \begin{align*}\sum_{n=0}^{p-1}\f{D_n}{(-8)^n}
 &=\sum_{n=0}^{p-1}D_n v\Ls{3+\sqrt{-6}}6^n\e
 -12\Big(x\cdot\f{3+\sqrt{-6}}6-\f{x+y}2\Big)^2\\&
  =\f 12(2x+y\sqrt {-6})^2\e 8x^2-2p-\f{p^2}{8x^2}\mod {p^3}.\end{align*}
\v2

\section*{9. Proofs of congruences involving $A_n$}
\par\q For $\tau\in H$ let
$$h(\tau)=\Ls{\eta(\tau)\eta(6\tau)}{\eta(2\tau)\eta(3\tau)}^{12}=
\f{f_1(2\tau)^{12}f_2(3\tau)^{12}}{64}
=\Ls{f_1(2\tau)}{f_1(6\tau)}^{12}
=\Ls{f_2(3\tau)}{f_2(\tau)}^{12}.\eqno{(9.1)}$$
Since $f_1(\tau)=f_2(-\f 1{\tau})$ we see that
$$h(\tau)=\Ls{f_1(2\tau)}{f_1(6\tau)}^{12}
=\Ls{f_2(-\f 1{2\tau})}{f_2(-\f 1{6\tau})}^{12}=h\Big(-\f 1{6\tau}\Big).\eqno{(9.2)}$$

 \par{\bf Lemma 9.1} {\it We have
 $$h\Ls{2+\sqrt{-2}}6=1\q\t{and}\q h\Ls{3+\sqrt{-3}}6=-1.$$}
 \par{\it Proof.} Since $f_1(\tau)=f_2(-\f 1{\tau})$ and $f_2(\tau)=\zeta_{24}^{-1}f_2(\tau+1)$ we see that
 $$f_1\Ls{2+\sqrt{-2}}3=f_2\Ls{-2+\sqrt{-2}}2=\zeta_{24}^{-1}f_2\Ls{\sqrt{-2}}2
 =\zeta_{24}^{-1}f_1(\sqrt{-2})$$
 and $f_2\ls{2+\sqrt{-2}}2=\zeta_{24}f_2\ls{\sqrt{-2}}2$. Thus,
 $$f_1\Ls{2+\sqrt{-2}}3f_2\Ls{2+\sqrt{-2}}2=f_1(\sqrt{-2})
 f_2\Ls{\sqrt{-2}}2=\sqrt 2.$$
 It then follows that
 $$h\Ls{2+\sqrt{-2}}6=\f 1{64}\Big(f_1\Ls{2+\sqrt{-2}}3f_2\Ls{2+\sqrt{-2}}2\Big)^{12}
 =\f{(\sqrt{2})^{12}}{64}=1.$$

 \par Now we turn to proving $h\ls{3+\sqrt{-3}}6=-1.$ By Lemma 3.2, $f(\sqrt{-3})^{24}=\f 1{t(\f{\sqrt{-3}}2)}=256$. Thus, $f(\sqrt{-3})^{12}=\pm 16$. Since
 $$f(\sqrt{-3})^{12}=e^{-2\pi i\sqrt{-3}/4}
 \prod_{n=1}^{\infty}\big(1+e^{2\pi i\sqrt{-3}(n-\f 12)}\big)^{12}
 =e^{2\pi \sqrt{3}/4}
 \prod_{n=1}^{\infty}\big(1+e^{-2\pi \sqrt{3}(n-\f 12)}\big)^{12}>0,$$
 we must have $f(\sqrt{-3})^{12}=16$. From [7, p.261], $j\sls{1+\sqrt{-3}}2=0$ and so $\gamma_2(\f{1+\sqrt{-3}}2)=0$. Since $\gamma_2(\tau)=\f{f_2(\tau)^{24}+16}{f_2(\tau)^8}$, we get
 $f_2(\f{1+\sqrt{-3}}2)^{24}=-16$ and so $f_2(\f{1+\sqrt{-3}}2)^{12}=\pm 4i$. Since
 \begin{align*}\f{f_2\sls{1+\sqrt{-3}}{2}^{12}}i&=\f 1i2^6e^{\pi i(1+\sqrt{-3})/2}
 \prod_{n=1}^{\infty}\big(1+e^{\pi i(1+\sqrt{-3})n}\big)^{12}
 \\&=2^6e^{-\sqrt 3\pi/2}\prod_{n=1}^{\infty}\big(1+(-1)^ne^{-\sqrt 3\pi n}\big)^{12}>0,\end{align*}
 we must have
 $f_2(\f{1+\sqrt{-3}}2)^{12}=4i$.
  From the properties of Weber's functions,
 \begin{align*}&f_1\Ls{3+\sqrt{-3}}3=f_1\Big(1-\f 1{\sqrt{-3}}\Big)=\zeta_{48}^{-1}f\Big(-\f 1{\sqrt{-3}}\Big)=\zeta_{48}^{-1}f(\sqrt{-3}),
 \\&f_2\Ls{3+\sqrt{-3}}2=\zeta_{24}f_2\Ls{1+\sqrt{-3}}2.
 \end{align*}
 Thus,
 \begin{align*}h\Ls{3+\sqrt{-3}}6&=\f 1{64}f_1\Ls{3+\sqrt{-3}}3^{12}f_2\Ls{3+\sqrt{-3}}2^{12}
 \\&=\f 1{64}\cdot\zeta_{48}^{-12}f(\sqrt{-3})^{12}\cdot
 \zeta_{24}^{12}f_2\Ls{1+\sqrt{-3}}2^{12}\\&=\f{i\cdot 16\cdot 4i}{64}=-1.
 \end{align*}
\par{\bf Remark 9.1} Lemma 9.1 was stated by Beukers in [2, p.30] without proof.
\v2
\par{\bf Lemma 9.2}
 {\it The function $h(\tau)$ is a Hauptmodul for $\Gamma_{0}(6)^{+}$ with a unique pole at the cusp~$[\frac{1}{2}]$ and satisfies that for $\tau$ such that $h(\tau)$ is near~$0$,
$$
\sum_{n=0}^{\infty}A_{n}h(\tau)^{n}=\frac{\eta(2\tau)^{7}\eta(3\tau)^{7}}{\eta(\tau)^{5}\eta(6\tau)^{5}}.
$$
The eta product $\frac{\eta(2\tau)^{7}\eta(3\tau)^{7}}{\eta(\tau)^{5}\eta(6\tau)^{5}}$ is a weight~2 modular form for $\Gamma_{0}(6)^{+}$ with a unique zero at the cusp~$[\frac{1}{2}]$.}
\v2
{\it Proof.} By (2.1), for any $\begin{pmatrix}
    a&b\\c&d
\end{pmatrix}\in\Gamma_{0}(6)$, one has that
\begin{align*}
\eta\Big(\frac{a\tau+b}{c\tau+d}\Big)^{12}
&=(-1)^{ab+cd(1-a^{2})-ca+3c_{0}(a-1)+\frac{3}{2}(a^{2}-1)r}
(c\tau+d)^{6}\eta(\tau)^{12},\\
\eta\Big(3\frac{a\tau+b}{c\tau+d}\Big)^{12}
&=\eta\Big(\frac{a(3\tau)+3b}{(c/3)(3\tau)+d}\Big)^{12}
\\&=(-1)^{3ab+\frac{1}{3}cd(1-a^{2})-\frac{1}{3}ca+3c_{0}(a-1)
+\frac{3}{2}(a^{2}-1)r}
(c\tau+d)^{6}\eta(3\tau)^{12},\\
\eta\Big(2\frac{a\tau+b}{c\tau+d}\Big)^{12}
&=\eta\Big(\frac{a(2\tau)+2b}{(c/2)(2\tau)+d}\Big)^{12}
\\&=(-1)^{2ab+\frac{1}{2}cd(1-a^{2})-\frac{1}{2}ca+3c_{0}(a-1)
+\frac{3}{2}(a^{2}-1)(r-1)}
(c\tau+d)^{6}\eta(2\tau)^{12},\\
\eta\Big(6\frac{a\tau+b}{c\tau+d}\Big)^{12}
&=\eta\Big(\frac{a(6\tau)+6b}{(c/6)(6\tau)+d}\Big)^{12}
\\&=(-1)^{6ab+\frac{1}{6}cd(1-a^{2})-\frac{1}{6}ca+3c_{0}(a-1)
+\frac{3}{2}(a^{2}-1)(r-1)}
(c\tau+d)^{6}\eta(6\tau)^{12},
\end{align*}
and thus deduces that
$$
h\Big(\frac{a\tau+b}{c\tau+d}\Big)=h(\tau).$$
Also,
$h(-\f 1{6\tau})=h(\tau)$ by (9.2).
Therefore, $h(\tau)=q+O(q^{2})$ is a Hauptmodul for $\Gamma_{0}(6)^{+}$ with a simple zero at the cusp~$[i\infty]$. Since $h(\tau)$ has no zeros in $H$ as an infinite product, and $\Gamma_{0}(6)^{+}$ has two cusps $[i\infty]$ and $[\frac{1}{2}]$, by the Residue Theorem for compact Riemann surfaces, one can tell that $h(\tau)$ has a unique pole at the cusp~$[\frac{1}{2}]$.
Again, following the transformations
\begin{align*}
\eta\Big(\frac{a\tau+b}{c\tau+d}\Big)^{5}
&=\zeta_{24}^{5(ab+cd(1-a^{2})-ca+3c_{0}(a-1)+\frac{3}{2}(a^{2}-1)r)}
(c\tau+d)^{\f 52}\eta(\tau)^{5},\\
\eta\Big(3\frac{a\tau+b}{c\tau+d}\Big)^{7}
&=\eta\Big(\frac{a(3\tau)+3b}{(c/3)(3\tau)+d}\Big)^{7}
\\&=\zeta_{24}^{7(3ab+\frac{1}{3}cd(1-a^{2})-\frac{1}{3}ca+3c_{0}(a-1)
+\frac{3}{2}(a^{2}-1)r)}
(c\tau+d)^{\f 72}\eta(3\tau)^{7},\\
\eta\Big(2\frac{a\tau+b}{c\tau+d}\Big)^{7}
&=\eta\Big(\frac{a(2\tau)+2b}{(c/2)(2\tau)+d}\Big)^{7}
\\&=\zeta_{24}^{7(2ab+\frac{1}{2}cd(1-a^{2})-\frac{1}{2}ca+3c_{0}(a-1)
+\frac{3}{2}(a^{2}-1)(r-1))}
(c\tau+d)^{\f 72}\eta(2\tau)^{7},\\
\eta\Big(6\frac{a\tau+b}{c\tau+d}\Big)^{5}
&=\eta\Big(\frac{a(6\tau)+6b}{(c/6)(6\tau)+d}\Big)^{5}
\\&=\zeta_{24}^{5(6ab+\frac{1}{6}cd(1-a^{2})-\frac{1}{6}ca+3c_{0}(a-1)
+\frac{3}{2}(a^{2}-1)(r-1))}
(c\tau+d)^{\f 52}\eta(6\tau)^{5},
\end{align*}
one has that the eta product $F(\tau)=\frac{\eta(2\tau)^{7}\eta(3\tau)^{7}}{\eta(\tau)^{5}\eta(6\tau)^{5}}$ satifies that
$$
F\Big(\frac{a\tau+b}{c\tau+d}\Big)=(c\tau+d)^{2}F(\tau).
$$
This together with the transformation identity
$$
\frac{\eta(-1/3\tau)^{7}\eta(-1/2\tau)^{7}}{\eta(-1/6\tau)^{5}\eta(-1/\tau)^{5}}
=-\tau^{2}\frac{\eta(2\tau)^{7}\eta(3\tau)^{7}}{\eta(\tau)^{5}\eta(6\tau)^{5}}
$$
justifies that the eta product $\frac{\eta(2\tau)^{7}\eta(3\tau)^{7}}{\eta(\tau)^{5}\eta(6\tau)^{5}}$ is a weight~2 modular form for $\Gamma_{0}(6)^{+}$ with a unique zero at the cusp~$[\frac{1}{2}]$. Finally, details about the relation between these two functions can be found in [5, Section 6.9].

\v2
\par{\bf Proof of Theorem 1.29.} Since $p=x^2+3y^2=(x+y\sqrt{-3})(x-y\sqrt{-3})$, we may choose the sign of $x$ so that $x+y\sqrt{-3}\not\in p\Bbb Z_p$. Let $h(\tau)$ be given by (9.1). By Lemma 9.2, one can see that both
$h(\tau)$ and $\sum_{n=0}^{\infty}A_nh(\tau)^{n}
$ satisfy the assumptions in Beukers' theorem.
  Set
 $$a_n=A_n,\q c=6y,\q d=x-3y,\q \alpha=\f{3+\sqrt{-3}}6\q\t{and}\q N=6.$$
 Then clearly
 $(c,d)=1,\ N\mid c,$ $c\alpha\bar{\alpha}, c(\alpha+\bar{\alpha})\in\Bbb Z$,
$p=(c\alpha+d)(c\bar{\alpha}+d)$ and $c\alpha+d\notin p\Bbb Z_p$.
    By Lemma 9.1, $h\sls{3+\sqrt{-3}}6=-1$.
 Hence, applying Beukers' theorem(i) and Lemma 2.3 we obtain
 \begin{align*}\sum_{n=0}^{p-1}(-1)^nA_n
 &=\sum_{n=0}^{p-1}A_nh\Ls{3+\sqrt{-3}}6^n\e \Big(6y\cdot \f{3+\sqrt{-3}}6+x-3y\Big)^2
 \\&=(x+y\sqrt{-3})^2\e 4x^2-2p-\f{p^2}{4x^2}\mod {p^3}.\end{align*}

\par{\bf Remark 9.2}
Let $p$ be a prime such that $p\e 1,3\mod 8$ and so $p=x^2+2y^2$. Assume that $x+y\sqrt{-2}\not\in \Bbb Z_p$. In [2, Corollary 1.20], Beukers proved that
$$\sum_{n=0}^{p-1}A_n\e (x+y\sqrt{-2})^2\mod {p^3}\eqno{(9.3)}$$
using the fact $h\ls{2+\sqrt{-2}}6=1$, whose proof can now be found in Lemma~9.1. Thanks to Beukers' result, one can also verify the following conjecture of the first author[21]:
$$\sum_{n=0}^{p-1}A_n\e 4x^2-2p-\f{p^2}{4x^2}\mod {p^3},\eqno{(9.4)}$$
by combining (9.3) with Lemma 2.3.
Finally, we remark that
Z.W. Sun[27] conjectured the congruence for  $\sum_{n=0}^{p-1}A_n$ modulo $p^2$, and C. Wang and Z.W. Sun[29] proved that $\sum_{n=0}^{p-1}A_n\e 4x^2-2p\mod {p^2}$ for any prime $p=x^2+2y^2
\e 1,3\mod 8$.

\end{document}